\documentclass[11pt]{article}
\usepackage{fullpage}
\usepackage{url}
\usepackage{tabls}
\usepackage{epsf}
\usepackage{appendix}
\usepackage{ragged2e}
\usepackage[top=1in, bottom=1in, left=1in, right=1in]{geometry}
\usepackage{enumitem}
\usepackage{amsmath, amssymb, amsthm}
\usepackage{algorithmic}
\usepackage[]{algorithm2e}
\usepackage{mathrsfs}
\usepackage{tikz}
\setlist[itemize]{leftmargin=*}
\usepackage{caption}
\usepackage{tabularx}
\usepackage{authblk}
\usepackage{hyperref}
\usepackage{subfig}
\usepackage{multirow,booktabs}
\usepackage{xcolor,colortbl}

\newcommand{\apolloR}{APOLLO3\textsuperscript{\textregistered}}

\title{Adaptive mesh refinement on Cartesian meshes applied to the mixed finite element discretization of the multigroup neutron diffusion equations}
\author[1]{Patrick Ciarlet, Jr. \thanks{patrick.ciarlet@ensta-paris.fr} }

\author[2]{Minh-Hieu Do \thanks{minh-hieu.do@cea.fr}}
\author[2]{Fran\c{c}ois Madiot \thanks{francois.madiot@cea.fr}}
\affil[1]{POEMS, CNRS, INRIA, ENSTA Paris, Institut Polytechnique de Paris, 91120 Palaiseau, France. }
\affil[2]{Universit\'e Paris-Saclay, CEA, Service d'\'Etudes des R\'eacteurs et de Math\'ematiques Appliqu\'ees, 91191, Gif-sur-Yvette, France.}

\date{ }
\usepackage{tikz}
\usepackage{smartdiagram}
\usetikzlibrary{arrows,shapes,positioning,shadows,trees,chains}

\newtheorem{theorem}{Theorem}

\newtheorem{remark}{Remark}

\numberwithin{theorem}{section}
\numberwithin{lemma}{section}
\numberwithin{remark}{section}
\numberwithin{definition}{section}
\numberwithin{proposition}{section}
\numberwithin{corollary}{section}
\numberwithin{assumption}{section}

\newcommand{\ds}{\displaystyle}
\newcommand{\pa}{\partial}
\def\eps{{\varepsilon}}
\newcommand{\N}{\mathbb{N}}
\newcommand{\R}{\mathbb{R}}
\newcommand{\dive}{{\rm div}\,} 
 
\newcommand{\keffh}{{k_{{\rm eff},h}}} 
\newcommand{\grad}{{\bf grad}\,}

\def\nvec{{\bf n}}
\def\pvec{{\bf p}}
\def\qvec{{\bf q}}
\def\vvec{{\bf v}}
\def\xvec{{\bf x}}

\def\Pvec{{\bf P}}
\def\uvec{{\bf u}}
\def\vvec{{\bf v}}
\def\Hvec{{\bf H}}

\def\Lvec{{\bf L}}
\def\Qvec{{\bf Q}}

\def\Hcalmg{{{\cal H}}}
\def\Xcal{{\cal X}}
\def\Xcalmg{{{\cal X}}}

\def\Ical{{\bf \mathcal I}}

\def\Ncal{{\bf \mathcal N}}
\def\Pcal{{\bf \mathcal P}}

\def\Tcal{{\bf \mathcal T}}

\def\Dbb{{\mathbb D}}

\def\Lud{{\underline L}}
\def\Ludvec{{\bf {\underline L}}}

\def\Vmg{{V}}
\def\Wud{{\underline W}}

\def\Hud{{\underline H}}
\def\Hudvec{{\bf {\underline H}}}

\def\T{{\mathbb T}}
\def\D{{\mathbb D}}

\def\Q{{\mathbb Q}}

\tikzset{
  basic/.style  = {draw, text width=2cm, drop shadow, font=\sffamily, rectangle},
  root/.style   = {basic, rounded corners=2pt, thin, align=center,
                   fill=orange!50!red},
  level 2/.style = {basic, rounded corners=6pt, thin,align=center, fill=red!50,
                   text width=8em},
  level 3/.style = {basic, thin, align=left, fill=yellow!100, text width=6.5em}
}
 \tikzstyle{decision} = [diamond, draw, fill=yellow!50!orange, 
     text width=5em, text badly centered, node distance=2cm, inner sep=0pt,minimum height=1em]
   \tikzstyle{decision2} = [diamond, aspect=2, draw, fill=yellow!50!orange, 
     text width=8em, text badly centered, node distance=2cm, inner sep=0pt,minimum height=1em]  
 \tikzstyle{block} = [rectangle, draw, fill=blue!50, 
     text width=7em, text centered, rounded corners, minimum height=2em]
 \tikzstyle{line} = [draw, -latex']
 \tikzstyle{cloud} = [draw, ellipse,fill=pink!50,text width=5em,text badly centered, node distance=3cm,
     minimum height=2em]
 \tikzstyle{cloud2} = [draw, ellipse,fill=white, node distance=3cm,
 minimum height=2em]    

 \tikzstyle{cloud3} = [draw, circle, node distance=0.1cm,
 minimum height=0.001em]   
  \tikzstyle{block2} = [rectangle, draw, fill=green!50, 
  text width=12em, text centered, rounded corners, minimum height=2em]
   \tikzstyle{block3} = [rectangle, draw, fill=red!50, 
   text width=2.5em, text centered, rounded corners, minimum height=2em]
\tikzstyle{block4} = [rectangle, draw, fill=blue!50, 
     text width=13em, text badly centered, rounded corners, minimum height=2em]
\tikzstyle{block5} = [rectangle, draw, fill=blue!50, 
     text width=7em, text badly centered, rounded corners, minimum height=2em]
\tikzstyle{block6} = [rectangle, draw, fill=cyan!50, 
text width=7em, text badly centered, rounded corners, minimum height=2em]

\begin{document}
%
%
%
%
%
%
\maketitle
\justify 
\parskip 6pt plus 1 pt minus 1 pt

\begin{abstract}
The multigroup neutron diffusion equations are often used to model the neutron density at the nuclear reactor core scale.
Classically, these equations can be recast in a mixed variational form.
This chapter presents an adaptive mesh refinement approach 
based on {\em a posteriori} estimators. We focus on refinement strategies on Cartesian meshes, since such structures are common for nuclear reactor core applications. 

\textit{This preprint corresponds to the Chapter 19 of volume 60 in AAMS,
Advances in Applied Mechanics (to appear).}
\end{abstract}

\section{Introduction} 


In neutronics, one is interested in modelling the neutron density inside a reactor core. 
The neutron flux density in the reactor core is determined by solving the transport equation which depends on seven variables: space (3), direction (2), energy or modulus of the velocity (1), and  time (1). 
Due to the high dimensionality of the problem, the numerical resolution of this equation faces some challenges in terms of computational cost. In practice, the neutron flux density can be modeled by the diffusion equations at the reactor core scale. 

Achieving an accurate and efficient resolution of the discretization is challenging. Adaptive Mesh Refinement (AMR) strategies have emerged in the late seventies to tackle this crucial aspect. The objective is to reach a target accuracy, on the solution or on a functional output, at an optimal computational cost.
 In fact, AMR strategies can be categorized: $h$-refinement (mesh subdivision)
 \cite{Verfurth1994}\,; $p$-refinement (local increase of the order of approximation)
 \cite{Babuska1981}\,; or $r$-refinement (moving mesh) 
  \cite{CHR1999}.
  The above strategies can be mixed together, such as $hp$-refinement \cite{Babuvska1994,Daniel2018} and $hr$-refinement~\cite{Lang2003}. 
 In neutronics, we are typically in the case of heterogeneous coefficients which may induce some singularities in the solution of the problems at hand, 
that is a loss of regularity of the solution due to the discontinuities in the data. For this reason, we
focus on the mesh subdivision strategy, which is the most suitable solution to address this kind of difficulty.

Adaptive mesh refinement is based on an essential ingredient: {\em a posteriori} estimates.
In neutronics, {\em a posteriori} error control is an ongoing research topic. AMR based on {\em a posteriori} estimators for the transport equation have been studied in~\cite{FOURNIER2011221,fuhrer1997posteriori,LATHOUWERS20111373,LATHOUWERS20111867,OWENS2017352}. In~\cite[Section 3.3]{FOURNIER2011221}, the authors address the issue of the regularity assumption used for the a posteriori estimate. In~\cite{LATHOUWERS20111373,LATHOUWERS20111867,OWENS2017352}, the a posteriori estimates are based on the Dual Weighted Residual method where the guaranteed estimator involves an exact adjoint solution. In~\cite{fuhrer1997posteriori}, they devise a reliable estimate that relies on the definition of a dual problem and highlight a lack of efficiency due to the absence of stability in this dual problem. Rigorous estimates that do not require excess regularity,
 together with an adaptive mesh refinement strategy, have been proposed for the resolution of the source problem on the transport equation~\cite{dahmen2020adaptive}. Following this work, theoretical aspects on the eigenvalue problem has been addressed in~\cite{dahmen2023accuracy}. In these papers, the authors devise a numerical strategy relying on accuracy controlled operator evaluation, such as in~\cite{dahmen2020adaptive} for the resolution of the source problem.\\At the reactor core scale, using simplified models is common in the nuclear industry. Precisely, the simplified models can be the neutron diffusion model or the simplified transport model. In~\cite{CiDM23}, we derive rigorous {\em a posteriori} error estimates for mixed finite
element discretizations of the Neutron Diffusion equations, and propose an adaptive mesh refinement strategy that preserves the Cartesian structure. A first application of this approach to the criticality problem was performed in{~\cite{DoCiMa21} albeit with subefficient estimators}. Regarding the industrial context and specifically the numerical simulations, our approach is part of the development of a mixed finite element solver called MINOS~\cite{baudron2007} in the \apolloR \, code~\cite{mosca2024apollo}.

The outline of this chapter is as follows. Section~\ref{sec-notations} introduces some notations. Section~\ref{sec:MG_model} describes the physical model, the variational formulation and its discretization. Section~\ref{sec:algorithmic_motivations} gives some algorithmic motivations related to the MINOS solver~\cite{baudron2007}.  Section~\ref{sec:MG_error-estimator} details {\em a posteriori} error estimates. Section~\ref{sec:num_AMR} includes the algorithm of our Adaptive Mesh Refinement (AMR) approach. Section~\ref{sec:numerical_illustration} illustrates numerically our AMR strategy.

\section{Notations}\label{sec-notations}
We choose the same notations as in \cite{CiDM23}.
Throughout the paper, $C$ is used to denote a generic positive constant which is independent of the mesh size, the mesh and the quantities/fields of interest. We also use the shorthand notation $A\lesssim B$ for the inequality $A\leq C B$, where $A$ and $B$ are two scalar quantities, and $C$ is a generic constant. \\ 
Vector-valued (resp. tensor-valued) function spaces are written in boldface character (resp. blackboard characters)\,; for the latter, the index \textit{sym} indicates symmetric fields. Given an open set ${\mathcal O}\subset\R^{{d}}$, ${d}=1,2,3$, we use the notation $(\cdot,\cdot)_{0,{\mathcal O}}$ (respectively $\|\cdot\|_{0,{\mathcal O}}$) for the $L^2({\mathcal O})$ and $\Lvec^2({\mathcal O})=(L^2({\mathcal O}))^{{d}}$ inner products (resp. norms). 
\\ 
When the boundary $\pa {\mathcal O}$ is Lipschitz, $\nvec$ denotes the unit outward normal vector field to $\pa {\mathcal O}$. Finally, it is assumed that the reader is familiar with vector-valued function spaces related to the diffusion equation, such as $\Hvec(\dive;{\mathcal O})$, $\Hvec_0(\dive;{\mathcal O})$ etc. \\
Let $G\in\N\setminus\{0,1\}$.
Given a function space $W$, we denote by $\Wud$ the product space $W^G$.
We extend the notation $(\cdot,\cdot)_{0,{\mathcal O}}$ (respectively $\|\cdot\|_{0,{\mathcal O}}$) to the $\underline{L^2}({\mathcal O})$ and $\underline{\Lvec^2}({\mathcal O})$ inner products (resp. norms).\\
Specifically, we let $\Omega$ be a bounded, connected and open subset of $\R^d$ {for $d=2,3$}, having a Lipschitz boundary which is piecewise smooth. We split $\Omega$ into $N$ open, connected, disjoint parts $\{\Omega_i\}_{1\le i\le N}$ with Lipschitz, piecewise smooth boundaries $\overline{\Omega}=\cup_{1\le i\le N}\overline{\Omega_i}$: the set $\{\Omega_i\}_{1\le i\le N}$ is called a partition of $\Omega$. For a field $v$ defined over $\Omega$, we shall use the notations $v_i=v_{|\Omega_i}$, for $1\le i\le N$.\\
Given the partition $\{\Omega_i\}_{1\le i\le N}$ of $\Omega$, we introduce a function space with piecewise regular elements:
\[
\begin{array}{rcl}
\Pcal W^{1,\infty}(\Omega)&=&\left\{D\in L^{\infty}(\Omega)\,|\,D_i\in W^{1,\infty}(\Omega_i),\,1\le i\le N\right\}.
\end{array}
\]
To measure $\psi\in\Pcal W^{1,\infty}(\Omega)$, we use the natural norm $\|\psi\|_{\Pcal W^{1,\infty}(\Omega)} =\max_{i=1,N}\|\psi_i\|_{W^{1,\infty}(\Omega_i)}$.

\section{The neutron diffusion model}
\label{sec:MG_model}
{The neutron flux density in the reactor core is determined by solving the transport equation which depends on seven variables: space (3), direction (2), energy or modulus of the velocity (1), and time (1). It physically states the balance between, on the one hand, the emission of neutrons by fission and, on the other hand, the absorption, scattering, and leakage (at the boundaries) of neutrons. The most common discretization of the energy variable is the {multigroup} approximation where the energy domain is divided into subintervals called energy groups. In practice, the neutron flux density is usually modeled by the {multigroup} diffusion equations~\cite[Chapter 7]{DuHa76} at the reactor core scale.}

{
In many practical situations, only a steady-state solution is relevant and it requires to solve a generalized eigenvalue problem.
In the companion paper{~\cite{DoCiMa21}}, we perform adaptive mesh refinement for the {multigroup} diffusion case on this {so-called} criticality problem.} We present in this section some theoretical results underlying this approach on a source problem. For $G\ge2$, we let $\Ical_G := \{1,...,G\}$ be the set of energy groups.

 {Let $\qvec:= (q_{x}^g)_{x=1,d}^{g=1,G}$, 
$\qvec^g=(q_{x}^g)_{x=1,d}\in\R^d$ for $1\le g\le G$, 
and $\dive\qvec=(\dive_{\xvec} \qvec^g)^{g=1,G}\in\R^G$}.

Let $\T_e$ be the even removal matrix. At (almost) every point in $\Omega$, it is a matrix of $\R^{G\times G}$:
\begin{equation*}
 \forall(g,g')\in\Ical_G\times\Ical_G,\quad (\T_e)_{g,g'}=\left\{
\begin{aligned}
& \Sigma_{r,0}^g:=\Sigma_t^g-\Sigma_{s,0}^{g\to g} \quad\text{ if } g=g',\\
& -\Sigma_{s,0}^{g'\rightarrow g} \quad\text{ if }g\neq g'{,}
\end{aligned}
\right.
\end{equation*}
where $\Sigma_{s,0}^{g'\rightarrow g}$ are the Legendre moments of order 0 of the macroscopic scattering cross sections from energy group $g'$ to energy group $g$ and the coefficient $\Sigma_t^g$
 is the macroscopic total cross section of energy group $g$~\cite[Part 2, Chapter IV, Section A, p.124-127]{DuHa76}. 
We emphasize that, in general, the matrix $\T_e$ is not symmetric. \\
We denote $\D$ the diffusion matrix. At (almost) every point in $\Omega$, it is a diagonal matrix of $\R^{G\times G}$:
\begin{equation*}
\forall g\in\Ical_G,\quad \D_{g,g}={ D^g}.
\end{equation*}
{where $D^g$ is the scalar-valued diffusion coefficient of energy group $g$.}
Mathematically speaking, the coefficients of the measurable fields of matrices $\T_{e}$ and $\D$ are supposed to be such that:
\begin{equation}\label{MG_Pos}
\left\{
\begin{array}{ll}
(0)&
\forall g,g'\in\Ical_G,\ ({D}^g, \Sigma_{r,0}^g,\Sigma_{s,0}^{g'\to g})\in \Pcal W^{1,\infty}(\Omega)\times\Pcal W^{1,\infty}(\Omega)\times L^\infty(\Omega)\,,\\
(i)&\exists\,({D})_*,({D})^*>0,\ \forall\,g\in\Ical_G,\ ({D})_*\le {D}^g\le ({D})^*\mbox{ a.e. in }\Omega\,,\\
(ii)&\exists\,(\Sigma_{r,0})_*,(\Sigma_{r,0})^*>0,\ \forall\,g\in\Ical_G,\ (\Sigma_{r,0})_*\le \Sigma_{r,0}^g\le (\Sigma_{r,0})^*\mbox{ a.e. in }\Omega\,,\\
(iii)&\exists {\epsilon\in(0,(G-1)^{-1})},\ \forall\,g,g'\in\Ical_G, g'\neq g,\ |\Sigma_{s,0}^{g\rightarrow g'}|\leq\epsilon \Sigma_{r,0}^g\mbox{ a.e. in }\Omega.
\end{array}
\right.
\end{equation}
As a consequence of (\ref{MG_Pos}): the matrix $\T_{e}$ is (almost everywhere) strictly diagonally dominant (so it is invertible)\,; the matrix $\D$ is also invertible (almost everywhere).
\begin{remark}\label{rem_strictly-diag-dominant}
Hypothesis \eqref{MG_Pos}$-(iii)$ models accurately the core of a pressurized water reactor and, in this case, $\epsilon$ is usually
a small fraction of $(G-1)^{-1}$. We refer to Appendix~\ref{sec:appendix} for an illustration. So on every row of $\T_e$, the off-diagonal entries are much smaller than the diagonal entries. In particular, the inverse of $\T_e$ is well-approximated by the inverse of its diagonal.
\end{remark}
\begin{remark}\label{rem_upscattering}
In the absence of up-scattering, i.e. in the situation where for all $1\leq g<g'\leq G$ $\Sigma_{s,0}^{g'\rightarrow g}=0$, the matrix $\T_e$ is lower triangular by (energy group) block.
\end{remark}
We consider the following multigroup neutron diffusion equations, with vanishing Dirichlet boundary condition{, also denoted zero-flux boundary condition}.
The criticality problem (in its primal form) writes 
\begin{equation}\label{eq:MG_diff_criticality_primal}
\left\{\begin{array}{l}
\mbox{Find $(\phi,k_{\text{eff}})\in (\Hud^1_0(\Omega)\setminus\{0\})\times\R$ such that }
\\
k_{\text{eff}}\text{ is an eigenvalue with maximal modulus, and}\cr\\
-\dive (\D\,\grad\phi )+ \T_e\,\phi= \frac{1}{k_{\text{eff}}} M_f\phi\mbox{ in }\Omega,
\end{array}\right.
\end{equation}
where $M_f\in\left(\R\right)^{G\times G}$ is the fission matrix defined for all $(g,g')\in \Ical_G\times\Ical_G$ by
\begin{align*}
(M_f)_{g,g'}=\chi^g(\nu\Sigma_f)^{g'},
\end{align*} 
where $\chi^g$ is the total spectrum of energy group $g$,  $\nu^{g'}$ is the average number of neutrons emitted by fission of energy group $g'$ and $\Sigma_f^{g'}$ is the macroscopic fission cross section of energy group $g'$.
The coefficients of the fission matrix are assumed to be such that
\begin{align}
\left\{
\begin{aligned}
&\forall g \in \Ical_g, (({\nu\Sigma_f})^g,\chi^g)\in L^\infty(\Omega)\times L^\infty(\Omega),\\
&\forall g \in \Ical_g,\, 0\leq ({\nu\Sigma_f})^g \text{ a.e. in }\Omega, \text{ and } (\underline{\nu\Sigma_f})\neq { 0} \text{ a.e. in }\Omega.
\end{aligned}
\right.
 \label{eq:fission_assumption}
\end{align}
Under the assumptions \eqref{MG_Pos} and \eqref{eq:fission_assumption} on the coefficients, the eigenvalue $k_{\text{eff}}$ with the greatest modulus is simple, real and strictly positive (see \cite[Chapter XXI]{DLvol6}). Furthermore, the associated eigenfunction $\phi$ is real and positive a.e. in $\Omega$, up to a multiplicative constant. In neutronics, $\phi$ and $k_{\text{eff}}$ are respectively called the neutron flux and the effective multiplication factor. The physical state of the core reactor
is characterized by the effective multiplication factor: if $k_{\text{eff}}=1$: the nuclear chain
reaction is self-sustaining; if $k_{\text{eff}} > 1$: the chain reaction is diverging; if $k_{\text{eff}}<1$: the chain reaction
vanishes.
\\
 Given a source term $S_f\in \Lud^2(\Omega)$, the source problem writes,
\begin{equation}\label{eq:MG_diff_primal}
\left\{\begin{array}{l}
\mbox{Find $\phi\in \Hud^1_0(\Omega)$ such that}\cr
-\dive (\D\,\grad\phi )+ \T_e\,\phi= S_f\mbox{ in }\Omega.
\end{array}\right.
\end{equation}
Under the assumptions \eqref{MG_Pos} on the coefficients, the primal problem \eqref{eq:MG_diff_primal} is well-posed, in the sense that for all $S_f\in \Lud^2(\Omega)$, there exists one and only one solution $\phi\in \Hud^1_0(\Omega)$ that solves \eqref{eq:MG_diff_primal}, with the bound $\|\phi\|_{1,\Omega}\lesssim\,\|S_f\|_{0,\Omega}$. Provided that the coefficients $D^g$, for $1\leq g\leq G$, are piecewise smooth, the solution {has extra smoothness} (see eg. Proposition 1 in \cite{CiJK17} in the monogroup case). Instead of imposing a Dirichlet boundary condition on $\pa\Omega$, one can consider a Neumann or Fourier boundary condition $\mu_F^g\phi^g+(D^g\,\grad\phi^g)\cdot\nvec=0$, with $\mu_F^g\ge0$ for all $1\leq g\leq G$. Results are similar. In neutronics, the homogeneous Neumann boundary condition is called reflective boundary condition, while the Fourier boundary condition with $\mu_F^g=1/2$ is denoted vacuum boundary condition. 
\subsection{Mixed variational formulation}
\label{sec:MG_variational_formulation_discretization}
Let us introduce the function spaces:
\[\begin{array}{rcl}
 \Hcalmg &=&  \left\{\,\xi=(\qvec,\psi)\in\Ludvec^2(\Omega)\times \Lud^2(\Omega)\right\}\,,\ \|\xi\|_{\Hcalmg}=\left(\|\qvec\|_{0,\Omega}^2\,+\,\|\psi\|_{0,\Omega}^2\right)^{1/2}\,; \cr
\Xcalmg &=& \left\{\,\xi=(\qvec,\psi)\in\Hudvec(\dive,\Omega)\times \Lud^2(\Omega)\right\}\,,\ \|\xi\|_{\Xcalmg}=\left(\|\qvec\|_{\Hudvec(\dive,\Omega)}^2\,+\,\|\psi\|_{0,\Omega}^2\right)^{1/2}\,.
\end{array}\]
From now on, we use the notations: $\zeta=(\pvec,\phi)$ and $\xi=(\qvec,\psi)$.

The solution $\phi$ to \eqref{eq:MG_diff_primal} belongs to $\Hud^1(\Omega)$, so if one lets $\pvec=-\Dbb\,\grad\phi\in\Ludvec^2(\Omega)$, the multigroup neutron diffusion source problem may be written as:
\begin{equation}\label{eq:MG_diff-mixed}
\left\{\begin{array}{l}
\mbox{Find $(\pvec,\phi)\in\Hudvec(\dive,\Omega)\times \Hud^1_0(\Omega)$ such that}\cr
-\D^{-1}\,\pvec\,-\,\grad\phi=0\mbox{ in }\Omega,\cr
\dive\pvec\,+\,\T_e\phi=S_{f}\mbox{ in }\Omega.
\end{array}\right.
\end{equation}
Solving the mixed problem \eqref{eq:MG_diff-mixed} is equivalent to solving \eqref{eq:MG_diff_primal}.
The variational formulation for the mixed problem associated to problem \eqref{eq:MG_diff-mixed} writes 
\begin{equation}\label{eq:MG_VF-3-PC}
\left\{\begin{array}{l}
\mbox{Find $(\pvec,\phi)\in\Hudvec(\dive,\Omega)\times \Lud^2(\Omega)$ such that}\cr
\forall(\qvec,\psi)\in\Hudvec(\dive,\Omega)\times \Lud^2(\Omega),\quad c((\pvec,\phi),(\qvec,\psi))={(S_{f},\psi)_{0,\Omega}},
\end{array}\right.
\end{equation}
where the bilinear form $c$ is defined as
\begin{equation}\label{eq:MG_bilin-form-c-PC}
c\ :\ ((\pvec,\phi),(\qvec,\psi)) \mapsto 
-(\D^{-1}\,\pvec,\qvec)_{0,\Omega}
+(\phi,\dive\qvec)_{0,\Omega}
+(\psi,\dive\pvec)_{0,\Omega}
+(\T_e\,\phi,\psi)_{0,\Omega}.
\end{equation}
Note that $\D^{-1}$ is a measurable field of diagonal matrices, which positive diagonal entries that are uniformly bounded away from $0$. Hence, the form $c$
is continuous, bilinear and symmetric on $\Hudvec(\dive,\Omega)\times \Lud^2(\Omega)$, and
one may prove that the mixed formulation (\ref{eq:MG_VF-3-PC}) is well-posed, see theorem 3.16 in \cite{Giret2018}. 
As a matter of fact, the result is obtained by proving an inf-sup condition in $\Xcalmg=\Hudvec(\dive,\Omega)\times \Lud^2(\Omega)$, which we recall here.
\begin{theorem}\label{th:MG_VF-1-PC}
	Let $\D$ and $\T_e$ satisfy (\ref{MG_Pos}). Then, the bilinear symmetric form $c$ fulfills an inf-sup condition:
	\begin{equation}\label{eq:MG_inf-sup}
	\exists\eta>0,\quad
	\inf_{\zeta\in\Xcalmg}\,
	\sup_{\xi\in\Xcalmg}
	\ds\frac{c(\zeta,\xi)}{\|\zeta\|_{\Xcalmg}\,\|\xi\|_{\Xcalmg}}\geq\eta.
	\end{equation}
\end{theorem}
Similarly, the multigroup neutron diffusion criticality problem may be written as:
\begin{equation}\label{eq:MG_diff-mixed-criticality}
\left\{\begin{array}{l}
\mbox{Find $(\pvec,\phi,k_{\text{eff}})\in\Hudvec(\dive,\Omega)\times (\Hud^1_0(\Omega)\setminus\{0\})\times\R$ such that }\\
k_{\text{eff}}\text{ is an eigenvalue with maximal modulus, and}\cr\\
-\D^{-1}\,\pvec\,-\,\grad\phi=0\mbox{ in }\Omega,\cr
\dive\pvec\,+\,\T_e\phi=\frac{1}{k_{\text{eff}}} M_f\phi\mbox{ in }\Omega.
\end{array}\right.
\end{equation}
The variational formulation for the mixed problem associated to problem \eqref{eq:MG_diff-mixed-criticality} writes 
\begin{equation}\label{eq:MG_VF-3-PC-criticality}
\left\{\begin{array}{l}
\mbox{Find $(\pvec,\phi,k_{\text{eff}})\in\Hudvec(\dive,\Omega)\times (\Lud^2(\Omega)\setminus\{0\})\times\R$ such that}\\
k_{\text{eff}}\text{ is an eigenvalue with maximal modulus, and}\cr\\
\forall(\qvec,\psi)\in\Hudvec(\dive,\Omega)\times \Lud^2(\Omega),\quad c((\pvec,\phi),(\qvec,\psi))=\frac{1}{k_{\text{eff}}}(M_f\,\phi,\psi)_{0,\Omega}.
\end{array}\right.
\end{equation}

\subsection{Discretization}\label{ss-sec-FEM-PC}
We consider conforming discretizations of (\ref{eq:MG_VF-3-PC}). Let $(\Tcal_h)_h$ be a family of { meshes, made for instance of simplices, or of rectangles ($d=2$), resp. cuboids ($d=3$)}, indexed by a parameter $h$ equal to the largest diameter of elements of a given {mesh}. We introduce discrete, finite-dimensional, spaces indexed by $h$ as follows:
\[ \Qvec_h\subset\Hudvec(\dive,\Omega),\mbox{ and }L_h\subset \Lud^2(\Omega). \]
The conforming mixed discretization of the variational formulation (\ref{eq:MG_VF-3-PC}) is then:
\begin{equation}\label{eq:MG_VF-1h-PC}
\left\{\begin{array}{l}
\mbox{Find $(\pvec_h,\phi_h)\in\Qvec_h\times L_h$ such that}\cr
\forall(\qvec_h,\psi_h)\in \Qvec_h\times {L}_h,\quad c((\pvec_h,\phi_h),(\qvec_h,\psi_h)) = {(S_{f},\psi_h)_{0,\Omega}}.
\end{array}\right.
\end{equation}
Similarly, 
the discrete counterpart of the criticality problem writes
\begin{align}\label{eq:discrete_varf_criticality_pb}
\left\{\begin{array}{l}
\mbox{Find $(\pvec_h,\phi_h,\keffh)\in\Qvec_h\times (L_h\setminus\{0\})\times \R$ such that }\\
 \keffh \text{ is an eigenvalue with maximal modulus, and } \cr\\
\forall(\qvec_h,\psi_h)\in\Qvec_h\times L_h,\quad c((\pvec_h,\phi_h),(\qvec_h,\psi_h))=\frac{1}{\keffh}(
M_f\phi_h,\psi_h
)_{0,\Omega}.
\end{array}\right.
\end{align}
We define:
\begin{equation}\nonumber
\Xcalmg_h = \left\{\,\xi_h=(\qvec_h,\psi_h)\in\Qvec_h\times L_h\right\}\,,\mbox{ endowed with }\|\cdot\|_{\Xcalmg}\,.
\end{equation}
\begin{remark} {At some point, the discrete spaces are considered locally, i.e.. restricted to one element of the mesh. So, one introduces the local spaces $\Qvec_h(K)$, $L_h(K)$, $\Xcalmg_h(K)$ for every $K\in \Tcal_h$}.
\end{remark}
Explicit {\em a priori} error estimates may be derived, see eg. \cite[Section 5.2]{Giret2018}. \\
In this paper, we focus on the Raviart-Thomas-N\'ed\'elec (RTN) Finite Element~\cite{Nedelec1980,RaTh77}. \\
For {\em simplicial meshes}, that is meshes made of simplices,  the description of the finite element spaces RTN${}_k$, where $k\ge0$ is the order of the discretization can be found e.g. in \cite{BoBF13}. Similarly for the scalar fields of $L_h$ with the same order. \\ 
{ For {\em rectangular {or Cartesian} meshes},} a description of the Raviart-Thomas-N\'ed\'elec (RTN) finite element spaces can be found for instance in Section 4.2 of \cite{JaCi13}.
Let us detail the definition in dimension 3.
Let $P(K)$ be the set of polynomials on $K\in\Tcal_h$. We define the following subspace of $P(K)$,
\begin{align*}
\Q_{l,m,p}(K)=\{p(y_1,y_2,y_3) \in P(K) | \quad  p(y_1,y_2,y_3)=\sum_{i_1,i_2,i_3=0}^{l,m,p}\alpha_{i_1,i_2,i_3}y_1^{i_1}y_2^{i_2}y_3^{i_3}, \alpha_{i_1,i_2,i_3}\in\R \},
\end{align*}
with $l,m,p\in \mathbb{N}$ and $\mathbb{Q}_{l}(K)=\mathbb{Q}_{l,l,l}(K)$. For $k\ge0$, the definition of the RTN element of order $k$ is
\begin{align*}
\mathrm{RTN}_k(K)=\ [ \underline{\Q_{k+1,k,k}(K)}\times {\bf 0}\times {\bf 0} ]\oplus[ {\bf 0}\times  \underline{\Q_{k,k+1,k}(K)} \times {\bf 0}]\oplus[ {\bf 0}\times {\bf 0} \times  \underline{\Q_{k,k,k+1}(K)} ].
\end{align*}
The definitions of the finite element spaces RTN${}_k$ are then 
{
\begin{align*}
&\Qvec_{h}= \{\qvec_h\in \underline{\Hvec(\dive,\Omega)} \ |\ \forall K\in\Tcal_h,\ \qvec_h{}_{|K} \in  \mathrm{RTN}_k(K) \},\\
&L_h= \{\psi_h\in \underline{L^2(\Omega)} \ |\ \forall K\in\Tcal_h,\ \psi_h{}_{|K} \in  \underline{\mathbb{Q}_{k}(K)}\}.
\end{align*}
}

\section{Algorithmic motivations}
\label{sec:algorithmic_motivations}
In this section, we describe the algorithm of the MINOS solver~\cite{baudron2007}. Recall that in many practical situations, we are interested in solving the criticality problem~\eqref{eq:MG_diff_criticality_primal} or (\ref{eq:MG_diff-mixed-criticality}). 
For this reason, we detail the algorithm associated to the resolution of the discrete variational formulation~\eqref{eq:discrete_varf_criticality_pb}. The algorithm for the discrete source problem~\eqref{eq:MG_VF-1h-PC} follows a similar pattern.

For $1\leq g\leq G$, we denote $\Pvec_h^g$ and $\Phi_h^g$ the vector of the degrees of freedom of $\pvec_h^g$ and $\phi_h^g$. Furthermore, the current vector $\Pvec_h^g$ is sorted by direction, $1\leq x\leq d$. We may reformulate problem~\eqref{eq:discrete_varf_criticality_pb} in the following matrix formalism:\\
\begin{align*}
\text{Find }\left(Z_h:=\begin{bmatrix}
\Pvec_h^1&
\Phi_h^1&
\dots&
\Pvec_h^G&
\Phi_h^G&
\end{bmatrix}^T\ne0,\keffh\right) \text{ such that }
\mathbb{M}_h
Z_h
= \frac{1}{\keffh} 
\mathbb{F}_h
Z_h
,
\end{align*}
where $\keffh$ is an eigenvalue with maximal modulus, $\mathbb{M}_h$ and $\mathbb{F}_h$ are (energy group) block matrices such that for all $g,g'\in\{1,\dots,G\}$
\begin{align*}
(\mathbb{M}_h)_{g,g'}:=
\left\{
\begin{aligned}
&\begin{pmatrix}
-\mathbb{A}^{g,g}_h & \mathbb{B}_h\\
(\mathbb{B}_h)^T & \mathbb{T}_h^{g,g}
\end{pmatrix} &\text{ if $g'=g$}\\
&\begin{pmatrix}
0 & 0\\
0 & \mathbb{T}_h^{g,g'}
\end{pmatrix} &\text{ otherwise}
\end{aligned}
\right.\qquad;
\qquad
(\mathbb{F}_h)_{g,g'}:=
\begin{pmatrix}
0 & 0 \\
0 & \mathbb{F}_{\phi,h}^{g,g'} \\
\end{pmatrix}.
\end{align*}
The matrices $\mathbb{A}^{g,g}_h$, $\mathbb{T}_h^{g,g'}$, $\mathbb{F}_{\phi,h}^{g,g'}$, $1\leq g,g'\leq G$, and $\mathbb{B}_h$ are respectively associated to the following bilinear forms
\begin{equation*}\label{eq:bilin-form-a-PC}
a_{g,g}:\left\{
\begin{array}{rcl}
\Hvec(\dive,\Omega)\times \Hvec(\dive,\Omega)&\rightarrow&\R\\
(\uvec,\vvec)&\mapsto&\ds((D^g)^{-1}\,\uvec,\vvec)_{0,\Omega}
\end{array}
\right.;
\end{equation*}
\begin{equation*}\label{eq:bilin-form-t-PC}
t_{g,g'}:\left\{
\begin{array}{rcl}
L^2(\Omega)\times L^2(\Omega)&\rightarrow&\R\\
(u,v)&\mapsto&((\T_e)_{g,g'}u,v)_{0,\Omega}
\end{array}
\right.;
\end{equation*}
\begin{equation*}\label{eq:bilin-form-f-PC}
f_{g,g'}:\left\{
\begin{array}{rcl}
L^2(\Omega)\times L^2(\Omega)&\rightarrow&\R\\
(u,v)&\mapsto&(\chi^g(\nu\Sigma_f)^{g'}u,v)_{0,\Omega}
\end{array}
\right.;
\end{equation*}
\begin{equation*}\label{eq:bilin-form-b-PC}
b:\left\{
\begin{array}{rcl}
\Hvec(\dive,\Omega)\times L^2(\Omega)&\rightarrow&\R\\
{(\vvec, u)}&\mapsto&(u,\dive\vvec)_{0,\Omega}
\end{array}
\right..
\end{equation*}

A standard iterative resolution strategy is the inverse power method. For each iteration, a classical strategy is to solve the linear systems with matrix $\mathbb{M}_h$ via two nested level of iterations, the so-called outer and inner iterations. The slight difference is that we do not solve \textit{exactly} the linear systems. Instead, we perform outer iterations as detailed below, which results in approximate solves. It is known  that the convergence of the inverse power method is governed by the ratio between the first and second eigenvalue. Our resolution algorithm with approximate solves behaves similarly. In order to mitigate this effect, we apply the classical Chebyshev acceleration~\cite{varga1962}. This acceleration is an extrapolation method based on Chebyshev polynomials.
In the rest of this section, we will omit the subscript $h$.
\subsection{Outer iteration}
 The outer iteration is similar to a Gauss Seidel iteration on the energy groups. It is defined as follows. Let $Z_{M}= \begin{bmatrix}
\Pvec_{M}^1&
\Phi_{M}^1&
\dots&
\Pvec_{M}^G&
\Phi_{M}^G&
\end{bmatrix}^T$ and $k_M$ be the approximation of $Z$ and $k$ at the outer iteration number $M$. For increasing $g=1,\cdots,G$, we solve the linear systems
\begin{align*}
(\mathbb{M})_{g,g} \begin{bmatrix}
\Pvec^g_{M+1} \\ \Phi^g_{M+1}
\end{bmatrix} = 
\begin{bmatrix}
S^g_{M,\pvec} \\
S^g_{M,\phi}
\end{bmatrix}
\end{align*}
iteratively, where
\begin{align*}
\begin{bmatrix}
S^g_{M,\pvec} \\
S^g_{M,\phi}
\end{bmatrix}
 =
\frac{1}{k_M}\sum_{g'=1}^G(\mathbb{F})_{g,g'}
\begin{bmatrix}
 \Pvec^{g'}_{M}\\
 \Phi^{g'}_{M}
 \end{bmatrix} -\sum_{g'=1}^{g-1}(\mathbb{M})_{g,g'}
 \begin{bmatrix}
 \Pvec^{g'}_{M+1} \\
 \Phi^{g'}_{M+1}
 \end{bmatrix} 
 -\sum_{g'=g+1}^{G}(\mathbb{M})_{g,g'}
 \begin{bmatrix}
 \Pvec^{g'}_{M} \\
 \Phi^{g'}_{M}
 \end{bmatrix}.
\end{align*}
 These systems may be rewritten as 
\begin{align}
&(\mathbb{A}^{g,g}+\mathbb{B}(\mathbb{T}^{g,g})^{-1}\mathbb{B}^T)\Pvec^g_{M+1} =\mathbb{B}(\mathbb{T}^{g,g})^{-1}S^g_{M,\phi} -S^g_{M,\pvec}, 
\label{eq:pvec_g}\\
& \Phi^g_{M+1}= (\mathbb{T}^{g,g})^{-1}(S^g_{M,\phi}-\mathbb{B}^T\Pvec^g_{M+1}).  \label{eq:phi_g}
\end{align}
The eigenvalue is updated as follows
\begin{align}
k_{M+1} = k_M\frac{(\mathbb{F}Z_{M+1})^T(\mathbb{F}Z_{M+1})}{(\mathbb{F}Z_M)^T(\mathbb{F}Z_{M+1})}.
\label{eq:update_eigenvalue}
\end{align}
We define the residual of the method at the $M^{\text{th}}$ iteration by
\begin{align*}
\epsilon_M=\frac{\|k_{M+1}^{-1}\mathbb{F}Z_{M+1}-k_{M}^{-1}\mathbb{F}Z_{M}\|_{\ell^\infty}}{\|k_{M+1}^{-1}\mathbb{F}Z_{M+1}\|_{\ell^1}},
\end{align*}
where the norms $\|\cdot\|_{{\ell^\infty}}$ and  $\|\cdot\|_{{\ell^1}}$ are defined for all $Z\in \R^{\Ncal}$, $\Ncal>0$, by $\|Z\|_{{\ell^\infty}}=\max_{1\leq i\leq \Ncal} |Z_i|$ and $\|Z\|_{{\ell^1}}=\sum_{i=1}^{ \Ncal} |Z_i|$. The convergence criteria is based on this residual.

In the absence of up-scattering (cf. Remark~\ref{rem_upscattering}), the matrix $\mathbb{M}_h$ is lower triangular by (energy group) block and the outer iteration corresponds
to an exact resolution of the inverse power iteration.
\subsection{Inner iteration}
We apply an iterative scheme to solve~\eqref{eq:pvec_g}: it is a Gauss-Seidel iteration over the directional blocks that is called alternative direction iteration. These iterations are called inner iterations. We consider now that the matrices $\mathbb{T}^{g,g}, \, 1\leq g\leq G,$ are inexpensive to invert, i.e. that $(\mathbb{T}^{g,g})^{-1}, \, 1\leq g\leq G,$ are available. 
\begin{remark}
In the case of a $RTN$ finite element on rectangles ($d = 2$) or cuboids ($d = 3$), with well suited
finite element degrees of freedom, the resulting matrices $\mathbb{T}^{g,g}, \, 1\leq g\leq G,$ can be diagonal.
\end{remark}
Let us assume here that $d=3$. The matrices $\mathbb{A}^{g,g'}, 1\leq g,g'\leq G$, and  $\mathbb{B}$ are (directional) block matrices
\begin{align*}
&\mathbb{A}^{g,g'} = ( (\mathbb{A}^{g,g'}_{x,x'})_{1\leq x,x'\leq d} ),\\
&\mathbb{B}=((\mathbb{B}_x)_{1\leq x\leq d})^T.
\end{align*}
\begin{remark}
In the case of Cartesian meshes, $\mathbb{A}^{g,g'}$ is a (directional) block diagonal matrix.
\end{remark}
We omit here the subscript $M$ and the superscript $g$. Let $\Pvec_{J}$ be the approximation of $\Pvec$ at the inner iteration number $J$. For increasing $x=1,\cdots,d$, we solve the linear systems
\begin{align*}
(\mathbb{A}^{g,g}_{x,x}+\mathbb{B}_x(\mathbb{T}^{g,g})^{-1}(\mathbb{B}_x)^T)\Pvec_{x,J} = S_{x,J},
\end{align*}
where
\begin{align*}
S_{x,J} = Q -\sum_{x'=1}^{x-1}(\mathbb{A}^{g,g}_{x,x'}+ \mathbb{B}_x(\mathbb{T}^{g,g})^{-1}\mathbb{B}_{x'})\Pvec_{x',J+1} -\sum_{x'=x+1}^d(\mathbb{A}^{g,g}_{x,x'}+\mathbb{B}_x(\mathbb{T}^{g,g})^{-1}\mathbb{B}_{x'})\Pvec_{x',J},
\end{align*}
and $Q$ is the right-hand side of~\eqref{eq:pvec_g}.
\subsection{Global description}
We regroup in Algorithm~\ref{algorithm:alg2} the global strategy implemented in the MINOS solver.

\begin{algorithm}[H]
\caption{{The iterative strategy implemented in the MINOS solver }}
\label{algorithm:alg2}
\While{not converged}{
 Compute the source term $Q_{g}$\\
  $\ds Q_{g}=\sum_{g'=1}^G(\mathbb{F})_{g,g'}\begin{bmatrix}
 \Pvec^{g'}\\
 \Phi^{g'}
 \end{bmatrix}$\\
 \For{$g=1$ \KwTo $G$ }{
  update $(\mathbb{M})_{g,g'}$ contribution from other energy groups, $g'\neq g$, \\
 $\ds Q_{g}=Q_{g}-\sum_{g'\neq g}(\mathbb{M})_{g,g'}\begin{bmatrix}
 \Pvec^{g'}\\
 \Phi^{g'}
 \end{bmatrix}$\\
\While{not converged}
     {
         \For{$x=1$ \KwTo $d$ }
             {
                Compute the source term $Q_{g,x}$\\
                $\ds Q_{g,x}=\mathbb{B}_x(\mathbb{T}^{g,g})^{-1}Q^g_{\phi}-\sum_{x'\neq x}{(\mathbb{A}^{g,g}_{x,x'}+\mathbb{B}_x(\mathbb{T}^{g,g})^{-1}\mathbb{B}_{x'})}\Pvec_{x'}^g-Q^g_{\pvec_x}$\\
                Solve $\ds (\mathbb{A}^{g,g}_{x,x}+\mathbb{B}_x(\mathbb{T}^{g,g})^{-1}(\mathbb{B}_x)^T)\Pvec_x^g = Q_{g,x}$        
             }
          $\Phi^g= (\mathbb{T}^{g,g})^{-1}(Q^g_{\phi}-B^T\Pvec^g)$
         }
     }
     update eigenvalue:~\eqref{eq:update_eigenvalue}
 }
\end{algorithm}

\section{{\em A posteriori} error estimators for a mixed finite element discretization}\label{sec:MG_error-estimator}
In order to obtain {\em a posteriori} estimates, we use the {so-called} reconstruction of the discrete solution  {$\zeta_h=(\pvec_h,\phi_h)$} to (\ref{eq:MG_VF-1h-PC}). In what follows, we denote by $\tilde{\zeta}_h:=\tilde{\zeta}_h(\zeta_h)$ a reconstruction, and by $\eta:=\eta(\tilde{\zeta}_h)$ an estimator. Classically, our aim is to devise {\em reliable} and {\em efficient} estimators for the reconstructed error $\zeta-\tilde{\zeta}_h$, meaning that:
\[ \begin{array}{ll}
\|\zeta-\tilde{\zeta}_h\| \le \mathtt{C}\,\eta & \mbox{(reliability)} \cr
\eta \le \mathtt{c}\,\|\zeta-\tilde{\zeta}_h\| & \mbox{(efficiency)} \end{array} \]
where $\mathtt{C}$ and $\mathtt{c}$ are generic constants, and $\|\cdot\|$ is some norm to measure the error. In what follows, we consider that
\[\Vmg=\Hud^1_0(\Omega), \mbox{the original space of solutions, see (\ref{eq:MG_diff_primal}),}\]
is the {\em default} space of scalar reconstructed fields. We also introduce the {\em broken spaces}
\[ H^1(\Tcal_h) = \{\psi\in L^2(\Omega) \ |\  \psi\in H^1(K) , \forall K \in \Tcal_h \},\quad  \Hvec(\dive;\Tcal_h) = \{\qvec\in \Lvec^2(\Omega) \ |\  \qvec\in \Hvec(\dive;K) , \forall K \in \Tcal_h \}. \] 
A first approach has been suggested in~\cite[Chapter 8]{Giret2018}. 
The reconstruction $\tilde{\zeta}_h=(\tilde{\pvec}_h,\tilde{\phi}_h)$ is defined as 
\begin{align*}
&\tilde{\pvec}_h= \pvec_h \in \Qvec_h\subset \Hudvec(\dive;\Omega),\quad
\tilde{\phi}_h \in \Vmg.
\end{align*}
\begin{remark}
{For other boundary conditions, i.e. for a Neumann or Fourier boundary condition, the {\em default} space $\Vmg$ of scalar reconstructed fields would be equal to $\Hud^1(\Omega)$}.
\end{remark}
\subsection{Reconstruction of the discrete solution}\label{sec:reconstructions}
In this section, we list some approaches to devise a reconstruction {of the discrete solution $(\pvec_h,\phi_h)$, here obtained with the RTN${}_k$ finite element discretization, for $k\ge0$ }. We refer to~\cite[Section 5.1]{CiDM23} for a more exhaustive introduction of reconstruction approaches for RTN finite element spaces.  \\
{For illustrative purposes, we consider 
Cartesian meshes. Observe that the results presented in this section can be extended to the case of 
simplicial meshes~\cite{Vohralik2010}.
{Let us introduce some further notations, given such a mesh $\Tcal_h$.
The set of facets of $\Tcal_h$ is denoted $\mathcal{F}_h$, and it is split as $\mathcal{F}_h = \mathcal{F}_h^i \cup \mathcal{F}_h^e$, with $\mathcal{F}_h^e$ (resp. $\mathcal{F}_h^i$) being the set of boundary facets (resp. interior facets).
We denote by $\mathbb{Q}_k(\mathcal{T}_h)$ the space of piecewise polynomials of maximal degree $k$ of each variable on each rectangle (or cuboid) $K\in \mathcal{T}_h$. We let $\mathcal{V}_h^k$ be the set of interpolation points (or nodes) where the degrees of freedom of the $V$-conforming Lagrange Finite Element space of order $k$ are defined. And, for a node $a\in\mathcal{V}_h^k$, we denote by $\mathcal{T}_a$ the set of rectangles or cuboids $K$ such that $a\in K$.}}

\subsubsection{Average and average,+ reconstructions}
\label{sec:average_reconstruction}
To post-process, we consider first the averaging operator of the neutron flux 
$\mathcal{I}_{av} : \underline{\mathbb{Q}_k(\mathcal{T}_h)} \to \underline{\mathbb{Q}_{k+1}(\mathcal{T}_h)} \cap V$ such that
\begin{equation*}
{ \forall \phi_h\in\underline{\mathbb{Q}_k(\mathcal{T}_h)},\ \forall a\in\mathcal{V}_h^{k+1}},\quad \mathcal{I}_{av}(\phi_h)(a)= \frac{1}{|\mathcal{T}_a|} \displaystyle \sum_{K\in \mathcal{T}_a} \phi_h{}_{|K}(a).
\end{equation*}
The {\em average reconstruction}
is
{
\begin{equation}
\tilde{\zeta}_{av,h}=(\pvec_h,\mathcal{I}_{av}(\phi_h)).
\label{eq:reconstruction_average}
\end{equation}
}
Similarly, we define another averaging operator 
$\mathcal{I}_{av,+} : \underline{\mathbb{Q}_k(\mathcal{T}_h)} \to \underline{\mathbb{Q}_{k+2}(\mathcal{T}_h)} \cap V$ such that
\begin{equation*}
{ \forall \phi_h\in\underline{\mathbb{Q}_k(\mathcal{T}_h)},\ \forall a\in\mathcal{V}_h^{k+2}},\quad \mathcal{I}_{av,+}(\phi_h)(a)= \frac{1}{|\mathcal{T}_a|} \displaystyle \sum_{K\in \mathcal{T}_a} \phi_h{}_{|K}(a).
\end{equation*}
The {\em average,+ reconstruction}
is
{
\begin{equation}
\tilde{\zeta}_{av,h}=(\pvec_h,\mathcal{I}_{av,+}(\phi_h)).
\label{eq:reconstruction_average_plus}
\end{equation}
}

\subsubsection{Post-processing reconstruction}
We present here the approach proposed in~\cite{Arbogast1995}, { valid for any $k\ge0$}.
{It is shown there that the discrete solution $\zeta_h$ to~\eqref{eq:MG_VF-1h-PC}, set in $\Xcal_h=\Qvec_h\times L_h$, is also equal to the first argument of the solution}
of a hybrid formulation, where the constraint on the continuity of the normal trace of $\pvec_h$ is relaxed. Let
\[ \Lambda_h=\left\{\lambda_h\in \underline{L^2(\mathcal{F}^i_h)}\ |\ \exists \qvec_h \in \Qvec_h,\, \lambda_h{}_{|F}=\qvec_h\cdot \nvec_{|F},\, \forall F\in \mathcal{F}^i_h \right\}, \] be the space of the Lagrange multipliers and let $\tilde{\Xcal}_h=\Pi_{K\in\Tcal_h}\Xcal_h(K)$ be the unconstrained approximation space with the RTN${}_k$ local finite element spaces. By definition, $\Xcal_h$ is a strict subset of $\tilde{\Xcal}_h$. \\
{The hybrid formulation is:}
\begin{equation}\label{eq:hybrid_formulation}
\left\{\begin{array}{l}
\mbox{Find $(\zeta_h,\lambda_h)\in \tilde{\Xcal}_h\times\Lambda_h$ such that}\cr
\ds\forall (\xi_h,\mu_h)\in \tilde{\Xcal}_h\times\Lambda_h,\quad c(\zeta_h,\xi_h)-\sum_{F\in \mathcal{F}^i_h}\int_F \lambda_h[\qvec_h\cdot \nvec]+\sum_{F\in \mathcal{F}^i_h}\int_F \mu_h[\pvec_h\cdot \nvec]={(S_{f},\psi_h)_{0,\Omega}}.
\end{array}\right.
\end{equation}

Let $\Pi_{M_h}:\tilde{\Xcal}_h\times\Lambda_h \to M_h$ be the projection such that, {given $(\zeta_h,\lambda_h)\in \tilde{\Xcal}_h\times\Lambda_h$, 
$\widehat{\phi}_h=\Pi_{M_h}(\zeta_h,\lambda_h)$ is governed by}
 \[ 
 \forall (\psi_h,\mu_h)\in L_h\times \Lambda_h,\quad ({\T_e }{\widehat{\phi}_h},\psi_h)_{0,\Omega} + \sum_{F\in \mathcal{F}^i_h}\int_F { \widehat{\phi}_h}\mu_h = ({\T_e }\phi_h,\psi_h)_{0,\Omega} + \sum_{F\in \mathcal{F}^i_h}\int_F \lambda_h\mu_h.
 \]
See footnote\footnote{
\label{nbdp_Mh}
In dimension 3, the space $M_h$ is defined as 
$\displaystyle
M_h=\Pi_{K\in\Tcal_h}M_h(K),
$
where for all $K\in \Tcal_h$,
\begin{align*}
M_h(K)=\underline{\Q_{k+2,k,k}(K)}\oplus \underline{\Q_{k,k+2,k}(K)} \oplus \underline{\Q_{k,k,k+2}(K)}
\end{align*}
See~\cite{Arbogast1995} for the definition of {\em ad hoc} finite-dimensional spaces $M_h$ for other families and types of elements.
} for a possible choice of space $M_h$. \\
In the field of neutronics, the post-processing is defined in~\cite{fevotte2019} by $ \mathcal{I}^2_{\text{RTN}} : \tilde{\Xcal}_h\times\Lambda_h \to \underline{\mathbb{Q}_{k+2}(\mathcal{T}_h)} \cap V$ such that
\begin{equation*}
{ \forall (\zeta_h,\lambda_h)\in\tilde{\Xcal}_h\times\Lambda_h,\ \forall a\in\mathcal{V}_h^{k+2}},\quad \mathcal{I}^2_{\text{RTN}}{(\zeta_h,\lambda_h)}(a)= \frac{1}{|\mathcal{T}_a|} \displaystyle \sum_{K\in \mathcal{T}_a} { (\Pi_{M_h}(\zeta_h,\lambda_h))}{}_{|K}(a).
\end{equation*}
The chosen {\em post-processing reconstruction} is finally
\begin{equation}
\tilde{\zeta}_{\text{RTN},h}=(\pvec_h, \mathcal{I}^2_{\text{RTN}}{(\zeta_h,\lambda_h)}).
\label{eq:reconstruction_RTN_post-processing}
\end{equation}
{
\subsection{{\em A posteriori} error estimate}
\label{sec:MG_error_estimates}
Let us define 
\begin{align*}
&d_S(\zeta,\xi) = { (\D^{-1}\,\pvec,\qvec)_{0,\Omega}+(\T_e\phi,\psi)_{0,\Omega}}\\
&d(\zeta,\xi)= d_S(\zeta,\xi){ +(\psi,\dive\pvec)_{0,\Omega}-(\phi,\dive\qvec)_{0,\Omega}}= c(\zeta,(-\mathbf{q},\psi)). 
\end{align*}
The definition is extended to piecewise smooth fields on $\Tcal_h$ by replacing $\displaystyle\int_{\Omega}$ by $\displaystyle \sum_{K\in \Tcal_h}\int_K$. \\
Extending the approach in~\cite{CiDM23}, we propose 
to use a weighted $\Hudvec(\dive;\Tcal_h)\times \Lud^2(\Omega)$ norm. 
}
One may define 
\[ \delta_{e,K}^{max} = \max_{g\in\Ical_G}\sup_K((\T_e)_{g,g}),\ \delta_{e,K}^{min} = \min_{g\in\Ical_G}\inf_K((\T_e)_{g,g})\,;\ D_{K}^{max} = \max_{g\in\Ical_G}\sup_K(D^g),\ D_{K}^{min} = \min_{g\in\Ical_G}\inf_K(D^g).\]
We define the following norm on $\Xcalmg$, for all $\zeta\in \Xcalmg$,
\begin{align*}
&\|\zeta\|_{S,MG}^2=\sum_{K\in \Tcal_h}\|\D^{-1/2} \pvec\|_{0,K}^2 + \|\delta_{e}^{1/2} \phi\|_{0,K}^2 
 +\sum_{K\in \Tcal_h} h_K^2(D_{K}^{min})^{-1}\| (\,\dive\pvec)\|_{0,K}^2 ,
\end{align*}
where $\delta_{e}$ the diagonal part of $\T_{e}$, and $\delta_{e}^{1/2}$ is its square root. By design, the diagonal entries of $\delta_{e}$ are all strictly positive: $(\phi,\psi)\mapsto(\delta_{e}\phi,\psi)_{0,K}$ defines an inner product in $\Lud^2(K)$, so one can measure $\phi_{|K}$ in the associated norm $\phi\mapsto((\delta_e\phi,\phi)_{0,K})^{1/2}$, and finally $((\delta_e\phi,\phi)_{0,K})^{1/2} = \|\delta_{e}^{1/2} \phi\|_{0,K}$.
\begin{remark}
By contrast, in general, $\T_e$ is not a symmetric matrix, so $(\phi,\psi)\mapsto(\T_e\phi,\psi)_{0,K}$ does not define an inner product. On the other hand, since $\T_e$ is strictly diagonal dominant with small off-diagonal entries (cf. remark~\ref{rem_strictly-diag-dominant}), it is expected that $(\delta_e\phi,\phi)_{0,K}$ is an accurate approximation of $(\T_e\phi,\phi)_{0,K}$.
\end{remark}

We introduce the following ${\Xcalmg}_{K}$-local norm, for all $\zeta\in \Xcalmg$,
\begin{equation}\label{eq:MG_def_+_norm}
|\zeta|_{+,K}=\sup_{\xi\in {\Xcalmg}_{K}, \|\xi\|_{S,MG}\leq 1} d(\zeta,\xi),
\end{equation}
with $N(K)=\{K'\in\Tcal_h\ |\ \text{dim}_{H}({\partial K'\cap \partial K})=d-1\}$, where dim$_H$ is the Hausdorff dimension, and
\begin{align*}
&{\Xcalmg}_{K}=\left\{\zeta=(\pvec,\phi)\in {\Xcalmg}, \text{ Supp}(\phi)\subset K, \text{ Supp}({\pvec})\subset N(K)\right\}.
\end{align*}

Observe that the norm $\|\cdot\|_{S,MG}$ measures elements of $\Xcalmg$ in ${\underline{\Hvec}({\dive},{\Tcal}_h)}\times \Lud^2(\Omega)$ norm. This corresponds precisely to the energy norm (cf.~\cite[Chapter 8]{Giret2018}). 

{
\label{theorem:MG_local_norm_guaranteed}
Let $\zeta$ and $\zeta_h$ be respectively the solution to~\eqref{eq:MG_VF-3-PC} and~\eqref{eq:MG_VF-1h-PC}. Let $\tilde{\zeta}_h=(\pvec_h,\tilde{\phi}_h)$ be a reconstruction of $\zeta_h$ in $\Qvec_h\times \Vmg$. For any $K\in \Tcal_h$, we define the residual estimators
\begin{align}\label{eq:MG_def_residual_estimator}
{\eta_{r,K}}=\|\delta_{e}^{-1/2} (S_f-\dive \pvec_h - \T_e\tilde{\phi}_h)\|_{0,K},
\end{align}
the flux estimator
\begin{align}\label{eq:MG_def_flux_estimator}
{\eta_{f,K}} = \|\D^{1/2}(\D^{-1}\pvec_h+\textnormal{\textbf{grad}} \tilde{\phi}_h)\|_{0,K}.
\end{align}
Then it stands for all $K\in\Tcal_h$,
\begin{align}
&|\zeta-\tilde{\zeta}_h|_{+,K}\leq \left(\eta^2_{r,K} +\sum_{K'\in N(K)}\eta^2_{f,K'}\right)^{1/2}.\label{eq:MG_local1}
\end{align}
}
In order to state the next estimate on the local efficiency of the {\em a posteriori} error estimators, we will assume that the entries of the matrices $\T_{e}$ and $\D$ are piecewise polynomials on $\Tcal_h$ and $S_f\in L_h$.
In addition, we suppose that $\tilde{\phi}_h$ is piecewise polynomial on $\Tcal_h$. 
The former assumption can be restrictive in the sense that it imposes a constraint on the model. On the other hand, the latter assumption only depends on the choice of the reconstruction.
{Let $K \in \Tcal_h$ and let ${\eta}_{r,K}$ and ${\eta}_{f,K}$
be the residual estimators respectively given by~\eqref{eq:MG_def_residual_estimator} and~\eqref{eq:MG_def_flux_estimator}.} We have the following estimate
\begin{align}
{\eta}_{r,K}&\leq \mathtt{c} {\left(\frac{\delta_{e,K}^{max}}{\delta_{e,K}^{min}}\right)^{1/2}} |\zeta-\tilde{\zeta}_h|_{+,K}
, \label{eq:MG_norm+_local_efficiency_eta_r}\\
{\eta}_{f,K}&\leq \mathtt{C}{\left(\frac{D_{K}^{max}}{D_{K}^{min}}\right)^{1/2}}|\zeta-\tilde{\zeta}_h|_{+,K}
,\label{eq:MG_norm+_local_efficiency_eta_F}
\end{align}
where {$\mathtt{c}$ and $\mathtt{C}$ are constants which depend} only on the polynomial degree of $S_f$, $\T_{e}$, $d$, and the shape-regularity parameter $\kappa_K$,
with {the $|\cdot|_{+,K}$ norm defined as in (\ref{eq:MG_def_+_norm})}. 

The proof of these results for the multi-group diffusion equation is an extension of the study for the diffusion mono-group equation proposed in~\cite[Section 5.2.3]{CiDM23}. For the sake of brevity, details of the numerical analysis for the multigroup diffusion equation will be given in a companion paper.

\section{Adaptive mesh refinement}
\label{sec:num_AMR}
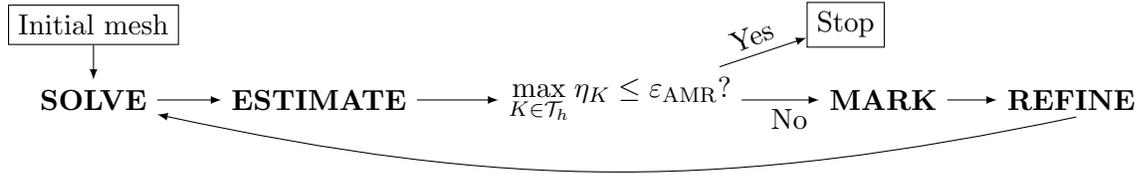
\begin{figure}[htbp]
\begin{tikzpicture}
\node [draw] (Im) at (0,2) {{Initial mesh}};
\node (S) at (0,1) {\textbf{SOLVE}};
\node (E) at (3,1) {\textbf{ESTIMATE}};
\node (T) at (7,1) {$\ds \max_{K \in\Tcal_{h}} \eta_K \leq \eps_{\text{AMR}}$?};
\node (M) at (10.5,1) {\textbf{MARK}};
\node (Mv) at (9.25,0.75) {No};
\node[rotate=21] (Mv) at (8.75,1.85) {Yes};
\node (R)  at (13,1) {\textbf{REFINE}};
\node (Rf)  at (13.2,0.8) { };
\node[draw] (Rd)  at (10,2) {Stop};
\draw[->,>=latex] (Im) -- (S);
\draw[->,>=latex] (S) -- (E);
\draw[->,>=latex] (E) -- (T);
\draw[->,>=latex] (T) -- (M);
\draw[->,>=latex] (M) -- (R);
\draw[->,>=latex] (T) -- (Rd) ;
\draw[->,>=latex] (Rf) to[bend left=12] (S);
\end{tikzpicture}
\caption{Description of the AMR process.}
    \label{fig:AMR_algo}
\end{figure}

In this Section, we recall a classical definition of an AMR strategy. 
This iterative process is divided into four modules as presented in Figure~\ref{fig:AMR_algo}, where $\eps_{\text{AMR}}>0$ is a user-defined parameter
%

\subsection{SOLVE module}
For the source problem, the  \textbf{SOLVE} module amounts to solving the discrete source problem~\eqref{eq:MG_VF-1h-PC}. In the case of the criticality problem, it corresponds to solving the discrete eigenvalue problem~\eqref{eq:discrete_varf_criticality_pb}.
\subsection{ESTIMATE module}
In the  \textbf{ESTIMATE} module, the  local error indicator $\eta_K$ is computed on each element $K\in \Tcal_h$. Using the {\em a posteriori} error estimate~\eqref{eq:MG_local1}, this error indicator is defined by
\begin{align*}
\eta_K:= \left(\eta^2_{r,K} +\sum_{K'\in N(K)}\eta^2_{f,K'}\right)^{1/2}.
\end{align*}
\subsection{MARK module}
The purpose of the  \textbf{MARK} module is to select a set of elements with large error: then, these elements are refined.
In other words, the marking strategy consists in selecting a set $S$ of elements of minimal cardinalF such that one has $$\eta(S) \simeq \theta \, \eta(\Tcal_{h}), \quad \text{where} \quad \eta(S) := \left(\sum_{K \in S} \eta_K^2\right)^{1/2}, \quad \mbox{resp.} \quad \eta(\Tcal_{h}) := \left(\sum_{K \in \Tcal_{h}} \eta_K^2\right)^{1/2}$$ and $\theta>0$ is a user-defined parameter.
According to~\cite[Section 6]{CiDM23}, an efficient 
strategy which preserves the Cartesian structure of the mesh is the \emph{direction} marker strategy. One selects for each direction $\mathbf{e}_x$, $x=1,\dots,d$, the smallest set of lines $L_{x}$ along that direction such that $\eta(L_{x}) \geq \theta \eta(\Tcal_{h})$. The resulting selected set is $\cup_{x=1,\dots,d} L_{x}$.
\subsection{REFINE module}
The \textbf{REFINE} module refines the mesh $\Tcal_{h}$ if the stopping criterion $\ds\max_{K \in\Tcal_{h}} \eta_K \leq \eps_{\text{AMR}}$ is not reached. 

\section{Numerical illustration}
\label{sec:numerical_illustration}


This section is dedicated to the numerical experiments where we apply the adaptive mesh refinement (AMR) described in Section~\ref{sec:num_AMR} to  industrial benchmark test cases.
 In fact, we compare the AMR process 
with respect to the choice of the flux reconstruction (involved in the definition of the total estimator):
\[ \mbox{(1) average; \quad (2) average,+; \quad (3) post-processing}. \] In particular, for the same threshold number $\eps_{\text{AMR}}>0$, the total number of elements denoted by $N_h$ is shown to illustrate the efficiency of the AMR strategy with different reconstruction methods.
In this section, we apply a RTN$_0$ finite element discretization and we consider the uniform Cartesian mesh with the mesh size $5 \text{cm}\times 5 \text{cm} \times 5 \text{cm}$ as the initial mesh for the AMR process, where the dimensions of $\Omega$ in the $x,y,z$-directions are in the order of $140 \text{cm}-320 \text{cm}$, see below. We slightly modify the boundary conditions of the benchmark test cases; we impose the zero-flux boundary condition instead of the vacuum boundary condition in our test cases.

\subsection{Takeda Model 2}

\begin{figure}[htbp]
\subfloat[]{
\includegraphics[]{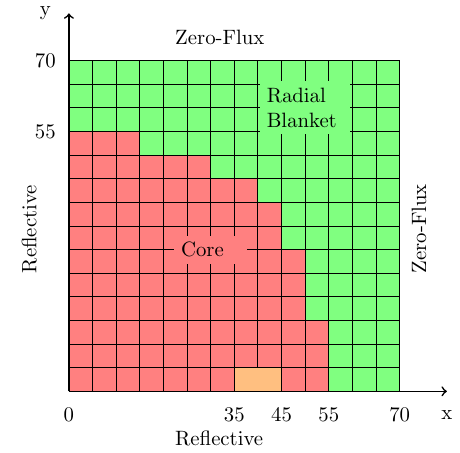}
}
\subfloat[]{
\includegraphics[]{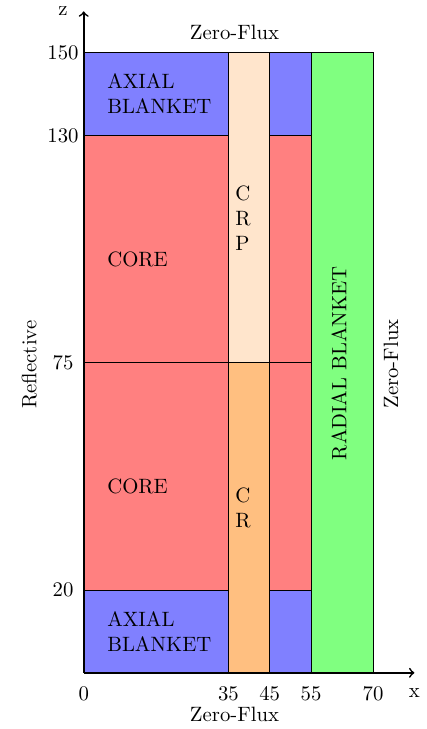}
}
\caption{Small FBR geometry}
\label{fig:geometry-model2}
\end{figure}

We aim to focus on the small Fast Breeder Reactor (FBR) (Model 2) of the 3D benchmark proposed  by Takeda et al. in \cite{takeda} with dimensions of $140 \text{cm}\times 140 \text{cm} \times 150 \text{cm}$  which includes a core region, radial and axial blankets, and a control rod region. The reactor core geometry (quarter reactor core) is shown in Figure \ref{fig:geometry-model2}.
The reactor core calculation is performed in 4 energy groups with the cross sections given in~\ref{sec:appendix}. In this problem, we would like to consider the following two cases:
\begin{itemize}
    \item Case 1: the control rod is withdrawn\,;
    \item Case 2: the control rod is half-inserted.
\end{itemize}

The so-called reference multiplication factors, computed on a uniform $2.5 \text{cm}\times 2.5 \text{cm} \times 2.5 \text{cm}$ mesh with RTN$_2$ elements are given in Table~\ref{tab:model2-keff}. According to Tables \ref{tab:model2-case1} and \ref{tab:model2-case2}, the AMR performs better if the estimators are computed with the average reconstruction than with the average,+ reconstruction.
This implies that using average operator with high order polynomial does not improve the accuracy of the estimator. On the other hand, the post-processing reconstruction reduces dramatically the total number of elements in the AMR process compared to the average reconstruction. Therefore, the post-processing reconstruction really improves the accuracy of the total estimator. Moreover, it is worth pointing out that the result depends strongly on the value of the parameter $\theta$ as shown in Tables \ref{tab:model2-case1-theta} and \ref{tab:model2-case2-theta}.
\\
Although the $\keffh$ obtained with the AMR processes of the average and average,+ reconstruction are closer to the reference multiplication factor, these processes are significanty less efficient with respect to the stopping criterion $\eps_{\text{AMR}}$. In addition, we observe that if one carries out more iterations of the AMR process of the post-processing reconstruction, one recovers a more accurate multiplication factor with fewer mesh elements. We refer to the values in italics in Table \ref{tab:model2-case1} where an error of 8 pcm (or $8.\, 10^{-5}$) is reached for $N_h\simeq 1.2 \, 10^{6}$, whereas $N_h> 2.3\, 10^7$ is needed for the other two methods to reach an error of 9 pcm (or $9. \,10^{-5}$) at best; similar observations are made with the values in italics in Table \ref{tab:model2-case2}.

\begin{figure}[htbp]
    \centering
   \subfloat[Case 1: $\min = 0.00011, \max = 0.08151$]{\includegraphics[width=0.45\linewidth]{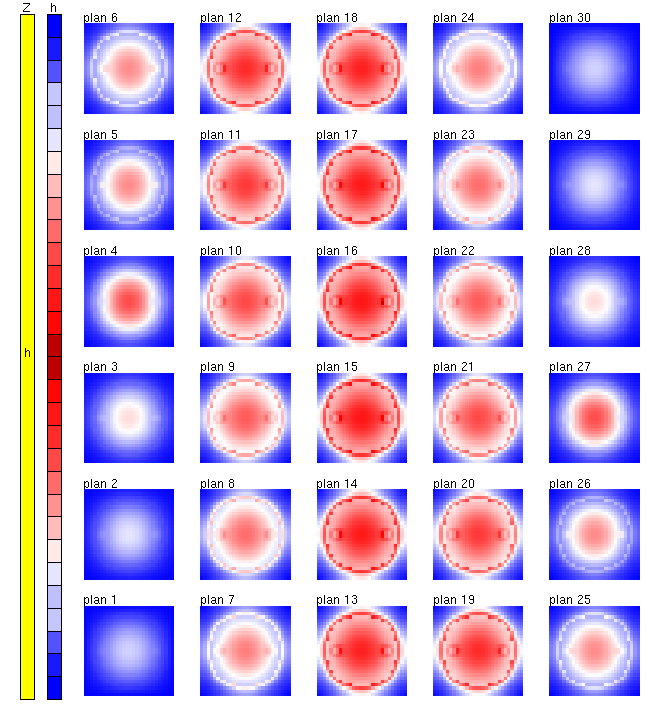}}\hspace{0.5cm}
   \subfloat[Case 2: $\min= 0.00011, \max = 0.45990$]{\includegraphics[width = 0.45\linewidth]{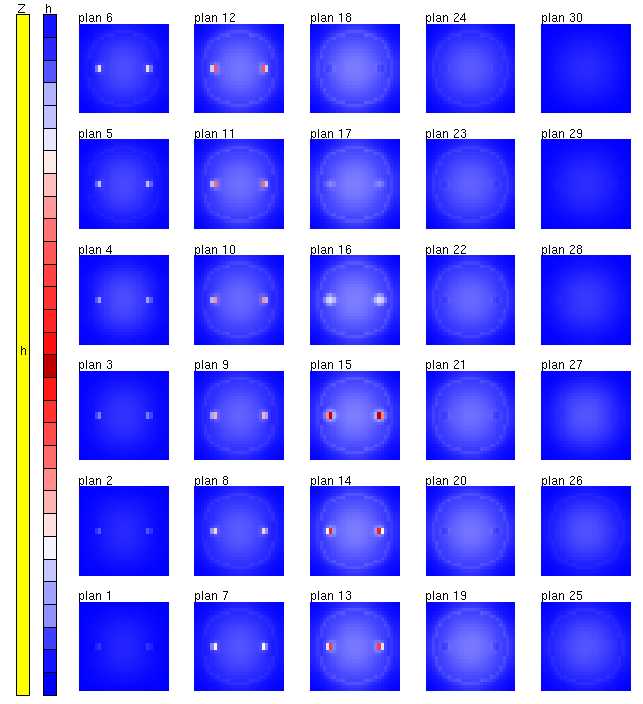}}
    \caption{Total estimator on the initial mesh for Model 2 - Post-processing reconstruction}
    \label{fig:model2-estimator}
\end{figure}

Figure \ref{fig:model2-estimator} illustrates the total estimator on the initial mesh for both cases. This figure clearly shows that the total estimator is large only at the rod positions for the test  with half inserted control rod (case 2) while it is large mainly on the interface between the Core and Radial Blanket for the test case where the control rod is fully withdrawn (case 1).

\begin{table}[htbp]
\centering
  	\begin{tabular}{c|cc}
			\toprule[\heavyrulewidth]\toprule[\heavyrulewidth]
			\textbf{Test case} & Case 1  & Case 2    \\
			\hline 
			$k_{\text{eff}}$ &  0.96361  & 0.94923  \\
			\hline 
		\end{tabular}
		\caption{The reference value of $k_{\text{eff}}$ for Model 2 }
		\label{tab:model2-keff}
\end{table}

\begin{table}[htbp]
  \centering 
\resizebox{\textwidth}{!}{  
  \begin{tabular}{cccccccccc} 
  \toprule[\heavyrulewidth]
  \toprule[\heavyrulewidth]
   & \multicolumn{3}{c}{Average} & \multicolumn{3}{c}{Average,+ } & \multicolumn{3}{c}{Post-processing} \\
   \cmidrule(r){2-4} \cmidrule(r){5-7} \cmidrule(r){8-10}
  \textbf{It} & \bf{$N_h$} &  $ \displaystyle \boldsymbol{\max_{K\in\Tcal_h}\eta_K }$ &  $\boldsymbol{\keffh}$ & \bf{$N_h$} &  $ \displaystyle \boldsymbol{\max_{K\in\Tcal_h}\eta_K }$ &  $\boldsymbol{\keffh}$& \bf{$N_h$}  &  $ \displaystyle \boldsymbol{\max_{K\in\Tcal_h}\eta_K }$ &  $\boldsymbol{\keffh}$ \\ 
  \midrule
0    &    23 520    &    0.40015    &    0.96405    &    23 520    &    0.49119    &    0.96405    &    23 520    &    0.08151    &    0.96405    \\
1    &    62 320    &    0.28764    &    0.96407    &    62 320    &    0.35711    &    0.96406    &    62 361    &    0.06078    &     0.96407    \\
2    &    165 996    &    0.23558    &    0.96410    &    157 092    &    0.35114    &    0.96410    &    188 328    &      \cellcolor{blue!25} 0.03626    &    0.96412    \\
3    &    473 550    &    0.13855    &    0.96413    &    450 000    &    0.19158    &    0.96413    &    \textit{564 570}    &    \textit{0.01318}    &    \textit{0.96397}    \\
4    &    921 600    &    0.07268    &    0.96374    &    902 400    &    0.11024   &    0.96374    &    \textit{1 168 128}    &    \textit{0.00993}    &    \textit{0.96369}    \\
5    &    2 315 163    &     \cellcolor{blue!25} 0.05644    &    0.96373    &    2 143 464    &    0.10117    &    0.96373    &       &      &      \\
6    &        &        &      &    6 289 630    &    \cellcolor{blue!25} 0.05894    &    0.96370    &      &        &        \\
  \bottomrule[\heavyrulewidth] 
  \end{tabular}
  }
    \caption{Takeda Model 2 case 1 with $\theta=0.5$ and $\eps_{\text{AMR}}=0.06$} 
  \label{tab:model2-case1}
\end{table}


\begin{table}[htbp]
  \centering 
\resizebox{\textwidth}{!}{  
  \begin{tabular}{cccccccccc} 
  \toprule[\heavyrulewidth]
  \toprule[\heavyrulewidth]
   & \multicolumn{3}{c}{Average} & \multicolumn{3}{c}{Average,+ } & \multicolumn{3}{c}{Post-processing} \\
   \cmidrule(r){2-4} \cmidrule(r){5-7} \cmidrule(r){8-10}
  \textbf{It} & \bf{$N_h$} &  $ \displaystyle \boldsymbol{\max_{K\in\Tcal_h}\eta_K }$ &  $\boldsymbol{\keffh}$ & \bf{$N_h$} &  $ \displaystyle \boldsymbol{\max_{K\in\Tcal_h}\eta_K }$ &  $\boldsymbol{\keffh}$& \bf{$N_h$}  &  $ \displaystyle \boldsymbol{\max_{K\in\Tcal_h}\eta_K }$ &  $\boldsymbol{\keffh}$ \\ 
  \midrule
0    &    23 520  &  1.67384    &    0.95016    &    23 520    &    1.10946    &    0.95016    &    23 520    &    0.45990    &    0.95016    \\
1    &    59 163 &   0.53247    &    0.94981    &    59 163    &    0.740026    &    0.94981    &    59 204    &    0.22517    &    0.95000    \\
2    &    159 848  &  0.24804    &    0.94976    &    148 257    &    0.31938   &    0.94972    &    175 230    &    0.06047   &    0.94983    \\
3    &    466 792  &  0.25122    &    0.94973    &    414 720    &    0.223137    &    0.94973    &    536 640    &    \cellcolor{blue!25} 0.02665    &    0.94974    \\
4    &    975 200  &  0.12797   &    0.94938    &    920 000    &    0.20069    &    0.94939    &    \textit{1 036 800}    &    \textit{0.01847}   &    \textit{0.94938}   \\
5    &    2 446 272  &  0.07737    &    0.94935    &    2 221 560    &    0.10953    &    0.94936    &    \textit{2 561 424}    &    \textit{0.01288}    &    \textit{0.94937}    \\
6    &    7 777 908  &  \cellcolor{blue!25}0.05029    &    0.94931    &    6 771 984    &    0.07823    &    0.94933    &    \textit{8 294 400}    &    \textit{0.00966}    &    \textit{0.94927}    \\
7    &      &      &       &    9 676 800    &   0.07822     &    0.94926    &        &        &   \\   8    &      &      &       &    21 546 096
    &   \cellcolor{blue!25} 0.040907
     &    0.94926    &        &        &        \\
  \bottomrule[\heavyrulewidth] 
  \end{tabular}
  }
    \caption{Takeda Model 2 case 2 with $\theta=0.5$ and $\eps_{\text{AMR}}=0.06$} 
  \label{tab:model2-case2}
\end{table}

\begin{table}[htbp]
  \centering 
  \begin{tabular}{cccccccccc} 
  \toprule[\heavyrulewidth]
  \toprule[\heavyrulewidth]
  & \multicolumn{3}{c}{$\theta=0.4$} & \multicolumn{3}{c}{$\theta=0.45$}& \multicolumn{3}{c}{$\theta=0.6$}  \\
  \cmidrule(r){2-4} \cmidrule(r){5-7} \cmidrule(r){8-10}
  \textbf{It} & \bf{$N_h$} &  $ \displaystyle \boldsymbol{\max_{K\in\Tcal_h}\eta_K }$ &  $\boldsymbol{\keffh}$ & \bf{$N_h$} &  $ \displaystyle \boldsymbol{\max_{K\in\Tcal_h}\eta_K }$ &  $\boldsymbol{\keffh}$ & \bf{$N_h$} &  $ \displaystyle \boldsymbol{\max_{K\in\Tcal_h}\eta_K }$ &  $\boldsymbol{\keffh}$\\ 
  \midrule
0    &    23 520    &    0.08151   &    0.96405  &    23 520    &    0.08151    &    0.96405   &    23 520    &    0.08151    &    0.96405    \\
1    &    53 391    &    0.06504    &    0.96406 &    57 760    &    0.06074    &    0.96406  &    75 768    &   \cellcolor{blue!25}  0.05246    &    0.96409    \\
2    &    130 000    &    \cellcolor{blue!25} 0.04628    &    0.96411 &    157 304    &     \cellcolor{blue!25} 0.04493    &    0.96411    &       &    &       \\
  \bottomrule[\heavyrulewidth] 
  \end{tabular}
    \caption{Takeda Model 2 case 1 with post-processing reconstruction and $\eps_{\text{AMR}}=0.06$} 
  \label{tab:model2-case1-theta}
\end{table}

\begin{table}[htbp]
  \centering 
  \begin{tabular}{cccccccccc} 
  \toprule[\heavyrulewidth]
  \toprule[\heavyrulewidth]
   & \multicolumn{3}{c}{$\theta=0.4$} &\multicolumn{3}{c}{$\theta=0.45$} & \multicolumn{3}{c}{$\theta=0.6$}  \\
   \cmidrule(r){2-4} \cmidrule(r){5-7} \cmidrule(r){8-10}
  \textbf{It} & \bf{$N_h$} &  $ \displaystyle \boldsymbol{\max_{K\in\Tcal_h}\eta_K }$ &  $\boldsymbol{\keffh}$ & \bf{$N_h$} &  $ \displaystyle \boldsymbol{\max_{K\in\Tcal_h}\eta_K }$ &  $\boldsymbol{\keffh}$& \bf{$N_h$} &  $ \displaystyle \boldsymbol{\max_{K\in\Tcal_h}\eta_K }$ &  $\boldsymbol{\keffh}$ \\ 
  \midrule
0    &    23 520    &    0.45990    &    0.95016  &    23 520    &    0.45990    &    0.95016   &    23 520    &    0.45990   &    0.95016    \\
1    &       50 544    &    0.24341    &    0.94999 &    54 760    &    0.22516    &    0.95000    & 72 283    &    0.20503    &    0.94994       \\
2        &    119 808    &    0.06190   &    0.94982 &    143 055    &    0.06189    &    0.94983  &    253 890    &   \cellcolor{blue!25} 0.05535    &    0.94984   \\
3   &    291 200    &   \cellcolor{blue!25}  0.05325    &    0.94983 &    383 116    &   \cellcolor{blue!25} 0.04903    &    0.94982   &       &      &         \\
  \bottomrule[\heavyrulewidth] 
  \end{tabular}
    \caption{Takeda Model 2 case 2 with post-processing reconstruction  and $\eps_{\text{AMR}}=0.06$} 
  \label{tab:model2-case2-theta}
\end{table}


\subsection{Takeda Model 3}
\begin{figure}[htbp]
    \centering
\subfloat[]{
\includegraphics[width=0.45\textwidth]{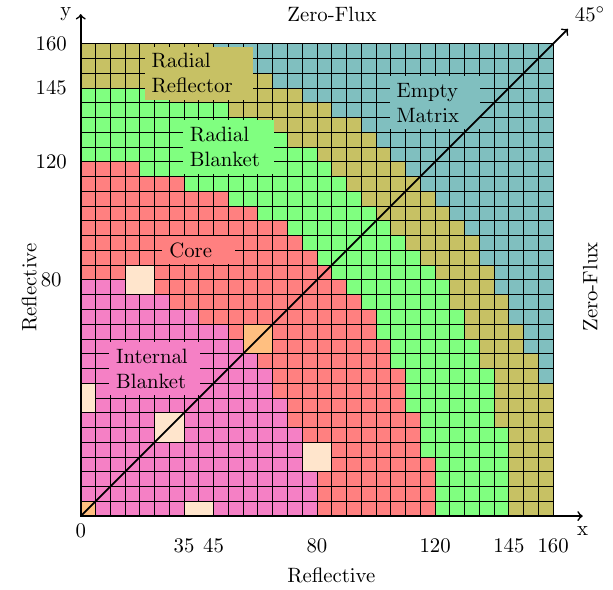}
    } \\
\subfloat[]{
\includegraphics[width=0.45\textwidth]{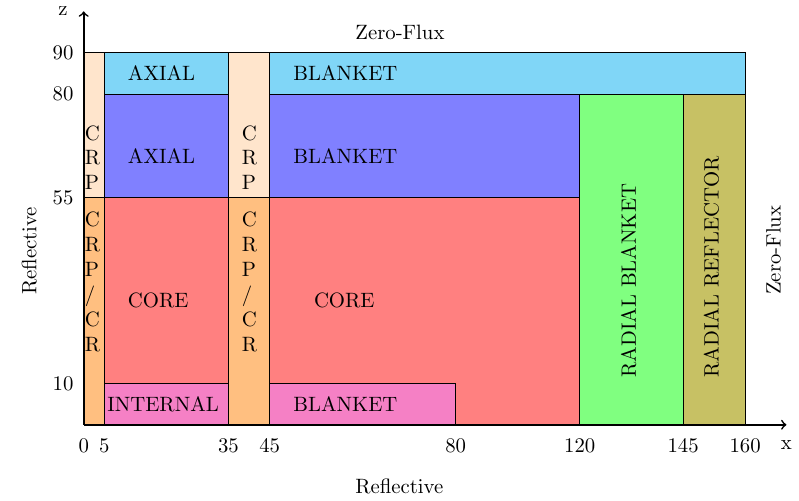}
}

    \caption{Configuration of Model 3 case 1 and case 2 }
    \label{fig:geometry-model3-case1}
\end{figure}

\begin{figure}[htbp]
    \centering
    \subfloat[]{
\includegraphics[width=0.45\textwidth]{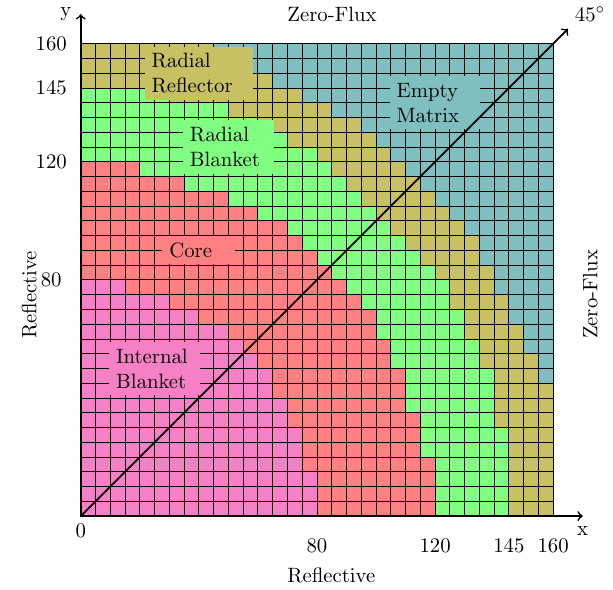}
    } \\
\subfloat[]{
\includegraphics[width=0.45\textwidth]{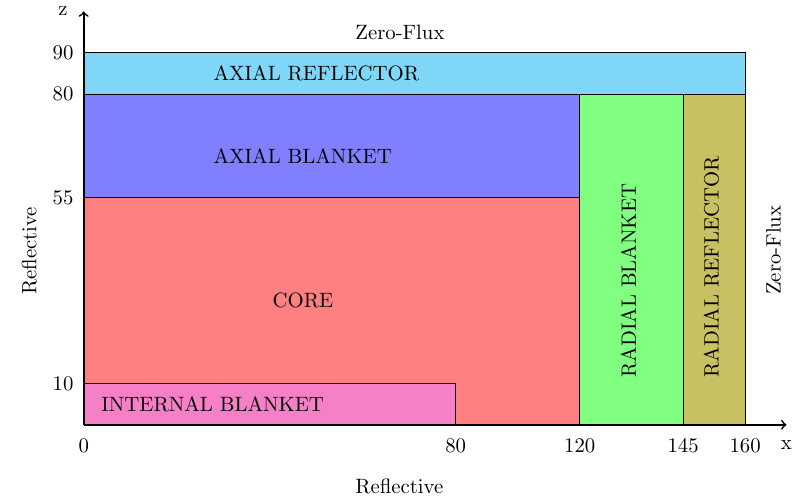}
}
    \caption{Configuration of Model 3  case 3}
    \label{fig:geometry-model3-case3}
\end{figure}

The third model of the benchmark proposed by Takeda et al. is an axially heterogeneous FBR with dimensions of $320 \text{cm}\times 320 \text{cm} \times 180 \text{cm}$. The geometry of the reactor core is shown in Figures \ref{fig:geometry-model3-case1} and \ref{fig:geometry-model3-case3}.
In particular, this reactor core is composed of internal blanket, core region, radial blanket, axial blanket, control rod region, radial reflector, axial reflector and empty matrix. 
The reactor core calculation is also performed in 4 groups with cross sections given in~\ref{sec:appendix} and the following three cases are considered:
\begin{itemize}
    \item Case 1: the control rods are inserted\,;
    \item Case 2: the control rods are withdrawn\,;
    \item Case 3: the control rods are replaced with core and/or blanket cells.
\end{itemize}

\begin{figure}[htbp]
    \centering
   \subfloat[Case 1: $\min = 3.28\times 10^{-6}, \max = 0.335723$ ]{\includegraphics[width=0.45\linewidth]{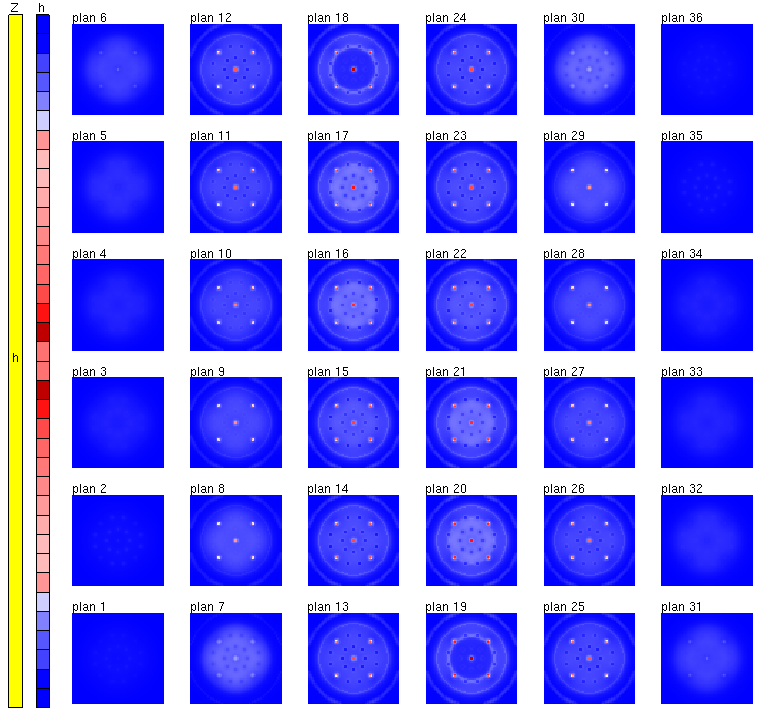}}\hspace{0.5cm}
   \subfloat[Case 2: $\min= 4.10\times 10^{-6}, \max = 0.07390$ ]{\includegraphics[width=0.45\linewidth]{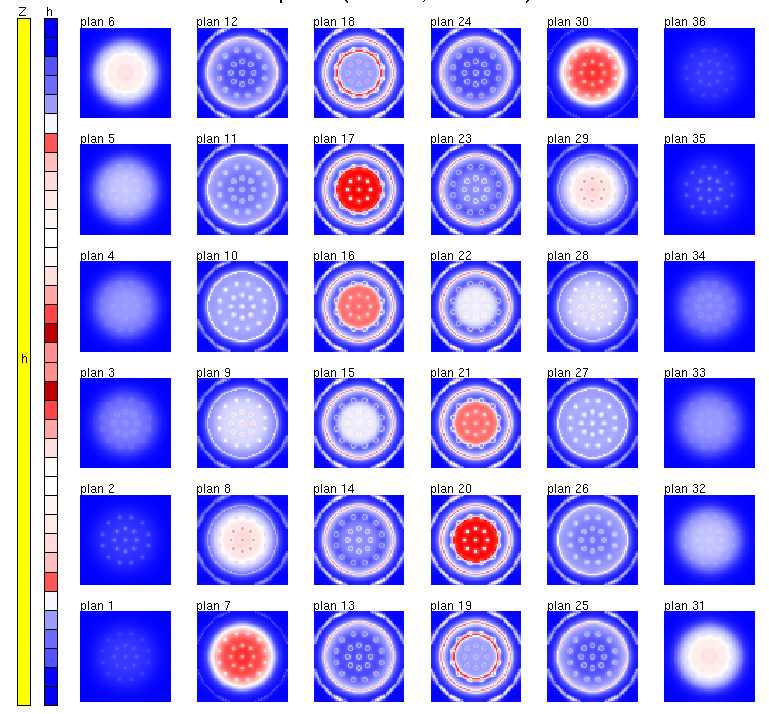}} \\
    \subfloat[Case 3: $\min =5.27\times 10^{-6}, \max=0.07585$]{\includegraphics[width=0.45\linewidth]{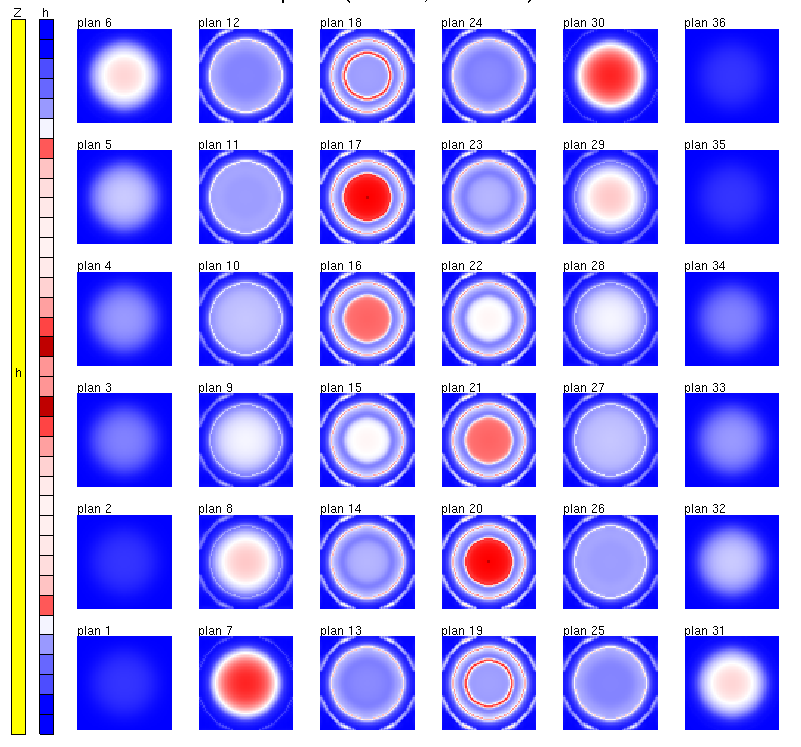}}
    \caption{Total estimator on the initial mesh for Model 3 - Post-processing reconstruction}
    \label{fig:model3-estimator}
\end{figure}

The so-called reference multiplication factors, computed on a uniform $2.5 \text{cm}\times 2.5 \text{cm} \times 2.5 \text{cm}$ mesh with RTN$_2$ elements are given in Table~\ref{tab:model3-keff}.
Since the average,+ reconstruction behaves similarly as discussed in the Takeda model 2, we only focus on the average and the post-processing reconstructions in this benchmark test case. Tables \ref{tab:model3-case1}, \ref{tab:model3-case2} and \ref{tab:model3-case3} show that the AMR process with the estimator computed with the post-processing reconstruction performs much better than with the average reconstruction.
Figure \ref{fig:model3-estimator} illustrates the behavior of the total estimator on the initial mesh for the model 3. As can be seen, it is large only at the control rods for the case 1, while it is large mainly at both the control rod positions and the interface among materials for case 2 and finally it is large only on the interface among materials for case 3. Moreover, Table \ref{tab:flux-model3} shows some initial refined meshes of the AMR process with the estimator computed with the post-processing reconstruction for all test cases, which are in good agreement with the behavior of the total estimator.

\begin{table}[htbp]
\centering
  	\begin{tabular}{c|ccc}
			\toprule[\heavyrulewidth]\toprule[\heavyrulewidth]
			\textbf{Test case} & Case 1  & Case 2  & Case 3  \\
			\hline 
			$k_{\text{eff}}$ &  0.96259  & 0.99333 & 1.01562  \\
			\hline 
		\end{tabular}
		\caption{The reference value of $k_{\text{eff}}$ for Model 3 }
		\label{tab:model3-keff}
\end{table}

\begin{table}[htbp]
  \centering 
  \begin{tabular}{ccccccc} 
  \toprule[\heavyrulewidth]
  \toprule[\heavyrulewidth]
   & \multicolumn{3}{c}{Average}  & \multicolumn{3}{c}{Post-processing} \\
   \cmidrule(r){2-4} \cmidrule(r){5-7} 
  \textbf{It} & \bf{$N_h$} &  $ \displaystyle \boldsymbol{\max_{K\in\Tcal_h}\eta_K }$ &  $\boldsymbol{\keffh}$ & \bf{$N_h$}  &  $ \displaystyle \boldsymbol{\max_{K\in\Tcal_h}\eta_K }$ &  $\boldsymbol{\keffh}$ \\ 
  \midrule
0    &    23 520    &    1.49448    &    0.96356                            &    23 520    &    0.33723    &    0.96356    \\
1    &    338 688    &    0.41143    &    0.96303                            &    354 025    &    0.08414    &    0.96307    \\
2    &    901 552    &    0.26117    &    0.96284                            &    962 948    &    \cellcolor{blue!25}0.05274    &    0.96290    \\
3    &    2 420 640    &    0.10432    &    0.96282                            &        &        &       \\
4    &    5 573 056    &    0.07457    &    0.96270                            &       &        &       \\
5    &    14 149 440    &    \cellcolor{blue!25}0.03725    &    0.96265                            &        &        &        \\
  \bottomrule[\heavyrulewidth] 
  \end{tabular}
  \caption{Takeda Model 3 case 1 with $\theta=0.5$ and $\eps_{\text{AMR}}=0.06$ } 
  \label{tab:model3-case1}
\end{table}



\begin{table}[htbp]
  \centering 
  \begin{tabular}{ccccccc} 
  \toprule[\heavyrulewidth]
  \toprule[\heavyrulewidth]
   & \multicolumn{3}{c}{Average} &  \multicolumn{3}{c}{Post-processing} \\
   \cmidrule(r){2-4} \cmidrule(r){5-7} 
  \textbf{It} & \bf{$N_h$} &  $ \displaystyle \boldsymbol{\max_{K\in\Tcal_h}\eta_K }$ &  $\boldsymbol{\keffh}$ & \bf{$N_h$}  &  $ \displaystyle \boldsymbol{\max_{K\in\Tcal_h}\eta_K }$ &  $\boldsymbol{\keffh}$ \\ 
  \midrule
0    &    23 520    &    0.26708    &    0.99348    &    23 520    &    0.07390    &    0.99348    \\
1    &    354 025    &    0.28711    &    0.99348    &    370 881    &    \cellcolor{blue!25}0.05321    &    0.99352    \\
2    &    948 787    &    0.17333    &    0.99347    &      &      &      \\
3    &    2 568 384    &    0.10356   &    0.99346    &        &      &       \\
4    &    5 537 792    &    \cellcolor{blue!25}0.04697    &    0.99336    &        &        &       \\
  \bottomrule[\heavyrulewidth] 
  \end{tabular}
    \caption{Takeda Model 3 case 2 with $\theta=0.5$ and $\eps_{\text{AMR}}=0.06$} 
  \label{tab:model3-case2}
\end{table}



\begin{table}[htbp]
  \centering 
  \begin{tabular}{cccccccccc} 
  \toprule[\heavyrulewidth]
  \toprule[\heavyrulewidth]
  & \multicolumn{3}{c}{Average} &  \multicolumn{3}{c}{Post-processing} \\
  \cmidrule(r){2-4} \cmidrule(r){5-7}
  \textbf{It} & \bf{$N_h$} &  $ \displaystyle \boldsymbol{\max_{K\in\Tcal_h}\eta_K }$ &  $\boldsymbol{\keffh}$ & \bf{$N_h$}  &  $ \displaystyle \boldsymbol{\max_{K\in\Tcal_h}\eta_K }$ &  $\boldsymbol{\keffh}$ \\ 
  \midrule
0	&	23 520	&	0.27265	&	1.01571	&	23 520	&	0.07585	&	1.01571	\\
1	&	354 025	&	0.18980	&	1.01573	&	37 0881	&	\cellcolor{blue!25}0.05058	&	1.01577	\\
2	&	934 626	&	0.16462	&	1.01574	&		&		&		\\
3	&	2 511 936	&	0.09597	&	1.01573	&		&		&		\\
4	&	5 537 792	&	\cellcolor{blue!25}0.04556	&	1.01564 	&		&		&		\\
  \bottomrule[\heavyrulewidth] 
  \end{tabular}
     \caption{Takeda Model 3 case 3 with $\theta=0.5$ and $\eps_{\text{AMR}}=0.06$ } 
  \label{tab:model3-case3}
\end{table}

\begin{table}[htbp]
\centering
\begin{tabular}{cccc}
&\textbf{Initial mesh} & \textbf{Iteration 1} & \textbf{Iteration 2} \\
\rotatebox{90}{\hspace{1.4cm}\small{\textbf{Case 1}}} 
&{\includegraphics[width=0.28\linewidth]{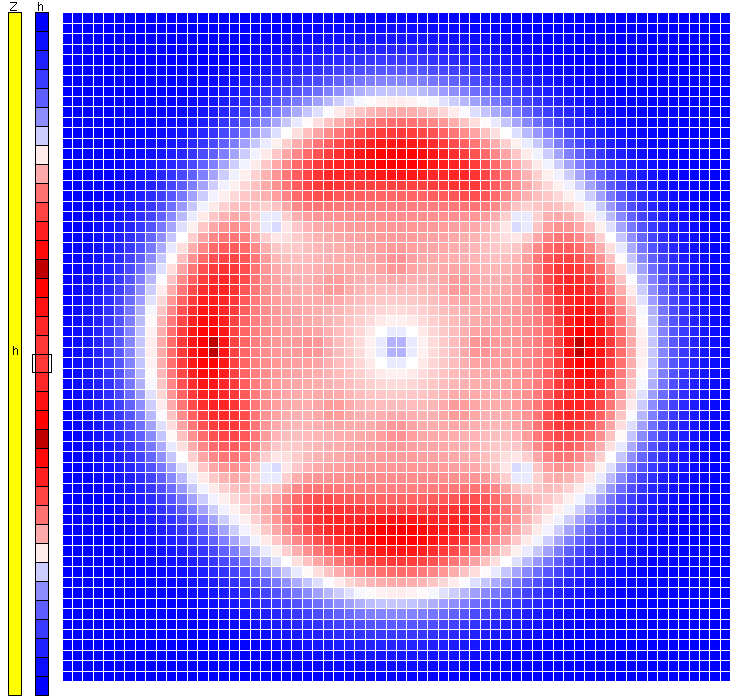}}&
{\includegraphics[width=0.285\linewidth]{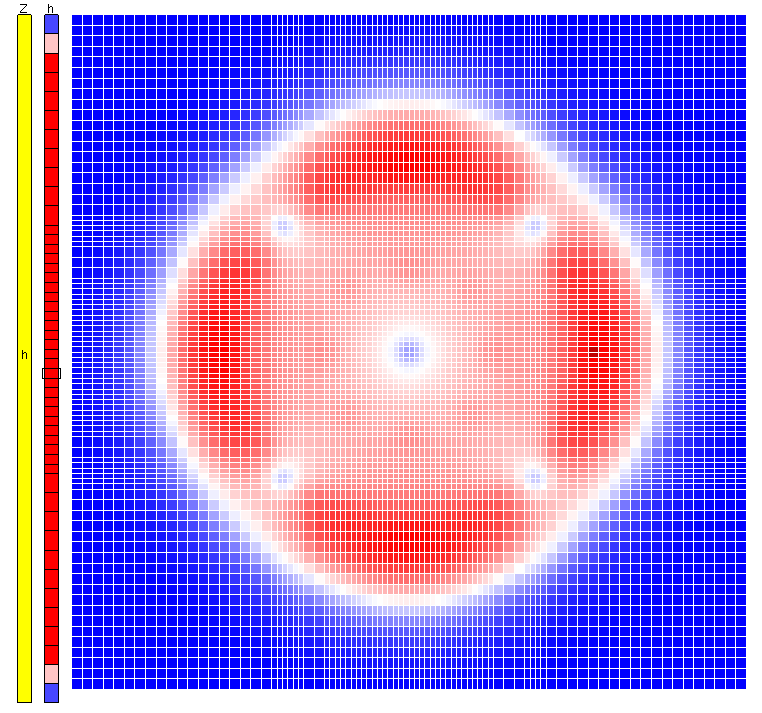}} &
{\includegraphics[width=0.28\linewidth]{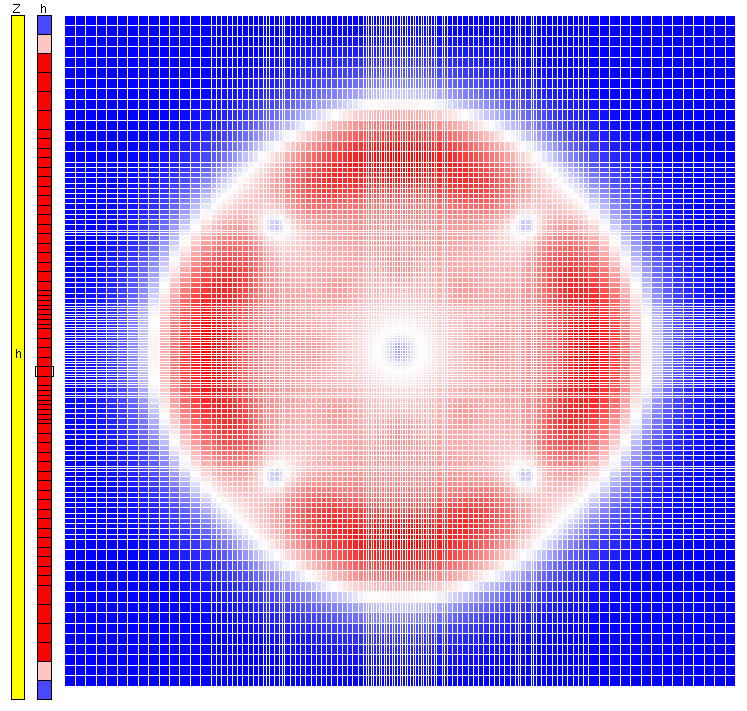}} \\
& \footnotesize{$\min= 2.755\, 10^{-5}, \max=1.379$}
 &\footnotesize{$\min=2.342\, 10^{-5}, \max=1.385$}
 & \footnotesize{$\min=2.343\, 10^{-5}, \max=1.390$}\\
 \rotatebox{90}{\hspace{1.4cm}\small{\textbf{Case 2}}}
 &{\includegraphics[width=0.28\linewidth]{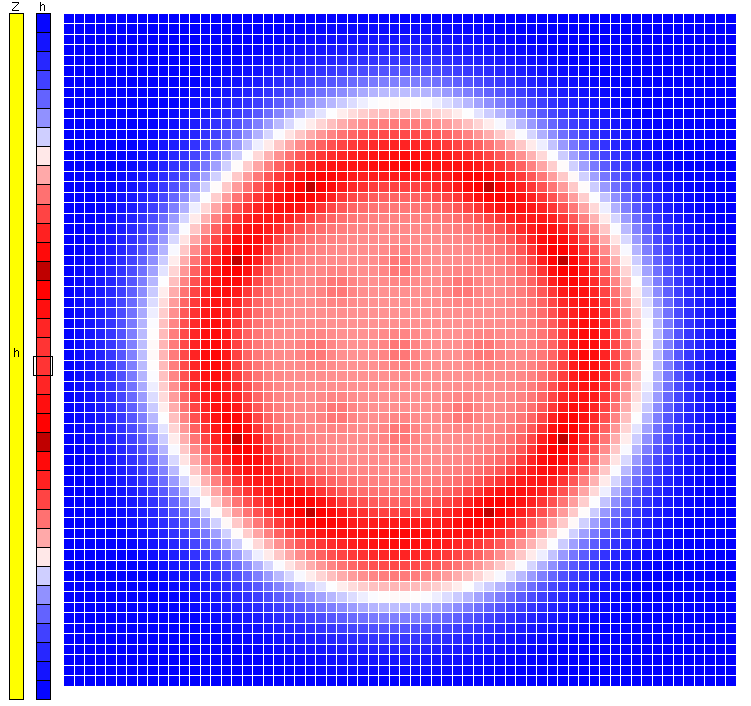}}&
{\includegraphics[width=0.277\linewidth]{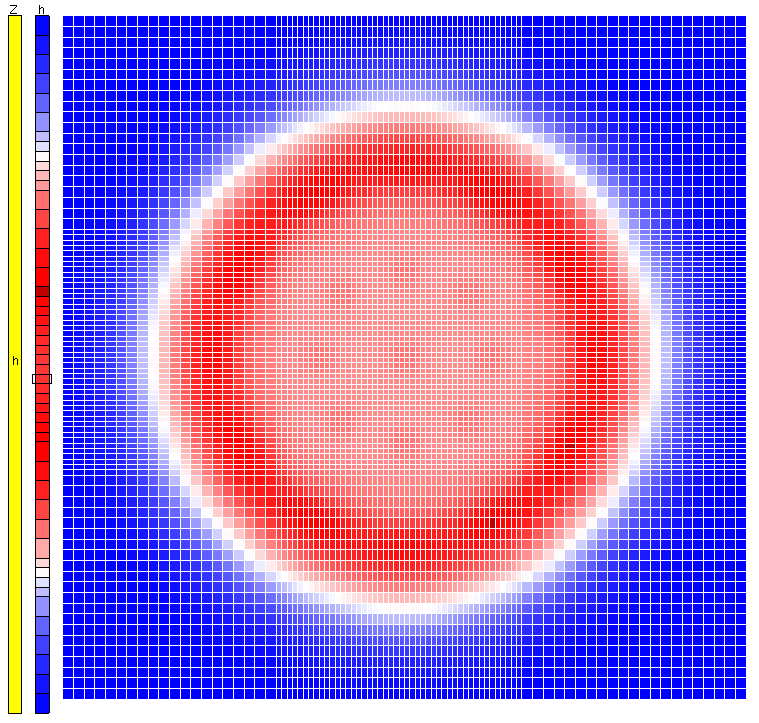}}&
{\includegraphics[width=0.28\linewidth]{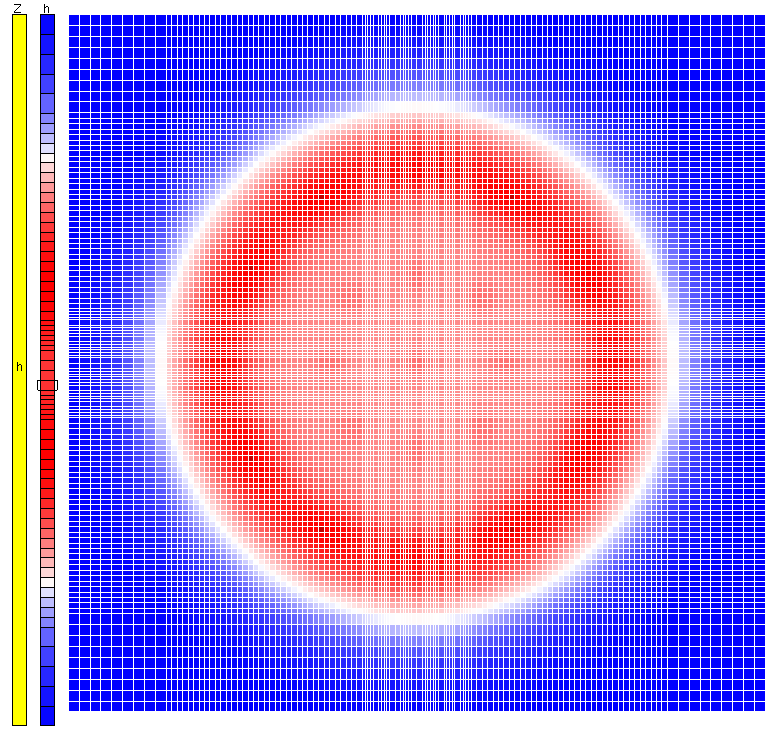}}\\
&\footnotesize{$\min=2.988\, 10^{-5}, \max=1.243$}
&\footnotesize{$\min=2.617\, 10^{-5}, \max=1.244$}
&\footnotesize{$\min=2.222\, 10^{-5}, \max=1.244$} \\
\rotatebox{90}{\hspace{1.4cm}\small{\textbf{Case 3}}}
&{\includegraphics[width=0.28\linewidth]{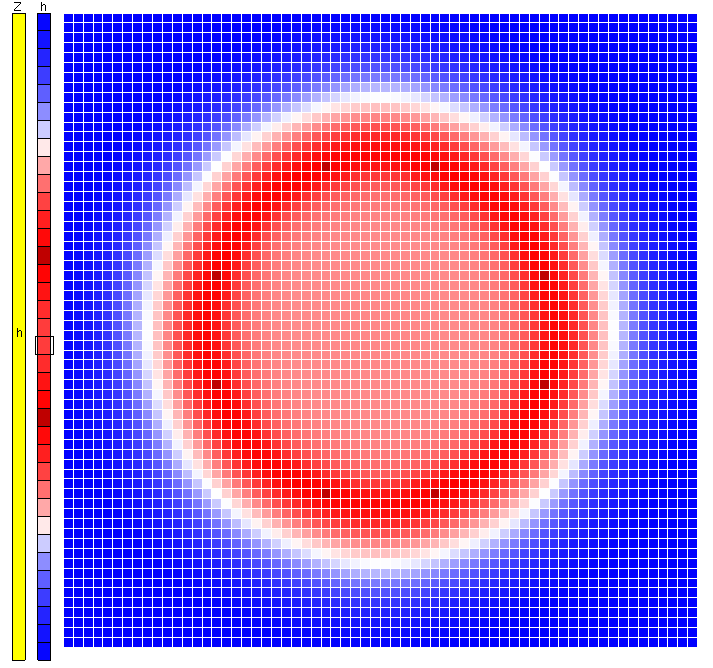}}&
{\includegraphics[width=0.28\linewidth]{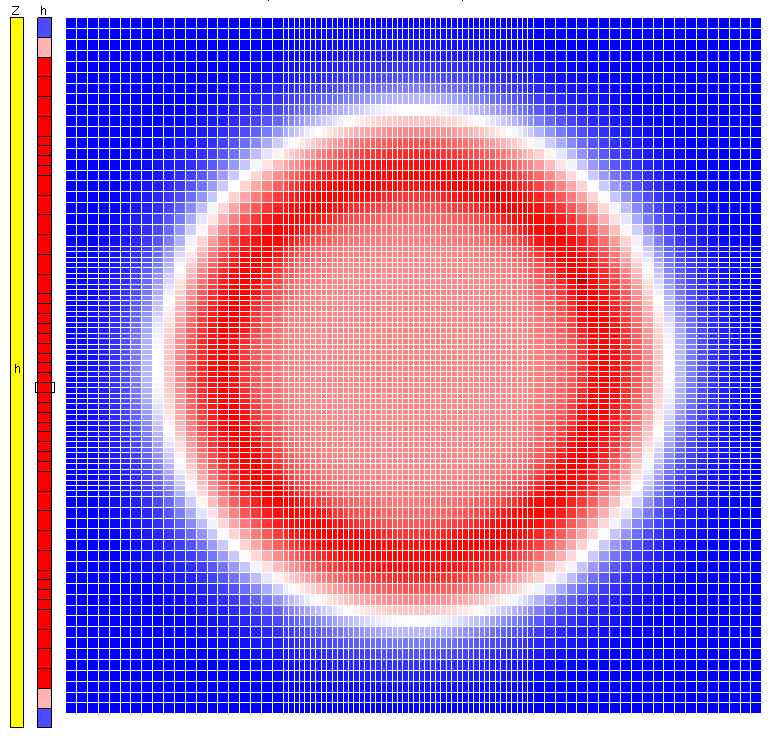}}&
{\includegraphics[width=0.28\linewidth]{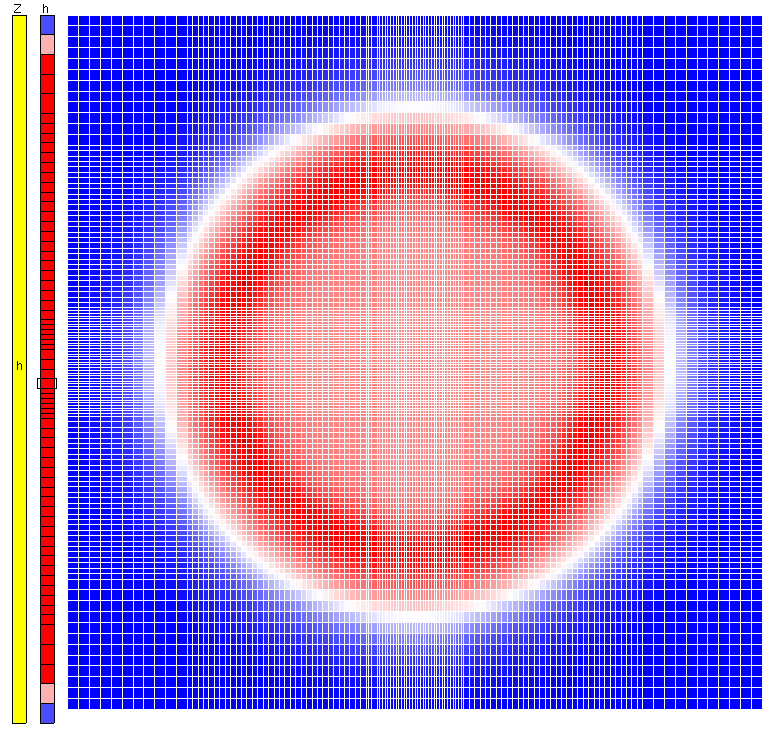}}\\
&\footnotesize{$\min=3.312\, 10^{-5}, \max=1.216$}
& \footnotesize{$\min=2.952\, 10^{-5}, \max=1.217$}
&\footnotesize{$\min=3.034\, 10^{-5}, \max=1.218$}\\
\end{tabular}
\caption{(Refined meshes) Flux of the second energy group of  model 3 at the central axial plane - Post-processing reconstruction}
\label{tab:flux-model3}
\end{table}
\newpage 

\section{Conclusion}
In this manuscript, we applied an AMR strategy to the mixed finite element discretization of the multigroup neutron diffusion equations. This approach is based on {\em a posteriori} estimators that are
both reliable and locally efficient. \\
We use a refinement strategy that preserves the Cartesian structure of the mesh. The potential of the approach is demonstrated on specific benchmarks coming from reactor physics application. We compare AMR processes with various a posteriori estimators. The post-processing reconstruction is shown to be more efficient than the average and average,+ reconstructions. We observe numerically that the AMR strategy is sensitive to the choice of the threshold parameter.\\
In a companion paper, we will consider a more general model, widely used for nuclear core simulations, the multigroup simplified tranport problem, for which we will provide a posteriori estimators. These estimators will be proven to be reliable and locally efficient.

\section*{Acknowledgements}
\apolloR \, is a registered trademark of CEA. We gratefully acknowledge EDF and Framatome for their long-term
partnership and their support.

\appendix
\section{Energy groups and cross-sections for the Takeda benchmark}
\label{sec:appendix}
For the sake of completeness, we detail in this appendix the energy groups and cross sections of the materials in the Takeda benchmarks.

\begin{table}[htbp]
  \centering 
  \begin{tabular}{cccc} 
  \toprule[\heavyrulewidth]\toprule[\heavyrulewidth]
  \multirow{2}{*}{Group} & \multicolumn{2}{c}{Energy range (eV)} & \multirow{2}{*}{Fission Spectrum $\chi$}\\
   & Upper & Lower &  \\
  \midrule
  1 & 1.00000$\, 10^7$ & 1.35340$\, 10^6$ & 0.583319 \\
  2 & 1.35340$\, 10^6$ & 8.65170$\, 10^{4}$ & 0.405450 \\
  3 & 8.65170$\, 10^{4}$ & 9.61120$\, 10^2$ & 0.011231 \\
  4 & 9.61120$\, 10^2$ & 1.00000$\, 10^{-5}$ & 0. \\
 
  \bottomrule[\heavyrulewidth] 
  \end{tabular}
    \caption{Energy range and Fission Spectrum} 
  \label{tab:data-fission-spectrum}
\end{table}

\begin{table}[htbp]
  \centering 
  \begin{tabular}{cccccc} 
  \toprule[\heavyrulewidth]
  \toprule[\heavyrulewidth]
   & \multicolumn{4}{c}{Core}  \\
   \cmidrule(r){2-5} 
  \textbf{Group} & \textbf{1} & \textbf{2} & \textbf{3} &  \textbf{4}  \\ 
  \midrule
  $\Sigma_t$  & 1.14568$\, 10^{-1}$ &  2.05177$\, 10^{-1}$ & 3.29381$\, 10^{-1}$ & 3.89810$\, 10^{-1}$\\
  $\nu\Sigma_f$ & 2.06063$\, 10^{-2}$ &  6.10571$\, 10^{-3}$ & 6.91403$\, 10^{-3}$ & 2.60689$\, 10^{-2}$ \\
  \midrule
  $\Sigma_{s,0}^{g\rightarrow 1}$ & 7.04326$\, 10^{-2}$ & 0. & 0. & 0.\\
  $\Sigma_{s,0}^{g\rightarrow 2}$ & 3.47967$\, 10^{-2}$ & 1.95443$\, 10^{-1}$ & 0. & 0.\\
  $\Sigma_{s,0}^{g\rightarrow 3}$ & 1.88282$\, 10^{-3}$ & 6.20863$\, 10^{-3}$ & 3.20586$\, 10^{-1}$ & 0.\\
  $\Sigma_{s,0}^{g\rightarrow 4}$ & 0. & 7.07208$\, 10^{-7}$ & 9.92975$\, 10^{-4}$ & 3.62360$\, 10^{-1}$\\
  
  \bottomrule[\heavyrulewidth] 
  \end{tabular}
    \caption{Cross sections in the Core region } 
  \label{tab:data-core}
\end{table}

\begin{table}[htbp]
  \centering 
  \begin{tabular}{cccccc} 
  \toprule[\heavyrulewidth]
  \toprule[\heavyrulewidth]
   & \multicolumn{4}{c}{Radial and Inner Blanket}  \\
   \cmidrule(r){2-5} 
  \textbf{Group} & \textbf{1} & \textbf{2} & \textbf{3} &  \textbf{4}  \\ 
  \midrule
  $\Sigma_t$        & 1.19648$\, 10^{-1}$ & 2.42195$\, 10^{-1}$ & 3.56476$\, 10^{-1}$ & 3.79433$\, 10^{-1}$\\
  $\nu\Sigma_f$     & 1.89496$\, 10^{-2}$ & 1.75265$\, 10^{-4}$ & 2.06978$\, 10^{-4}$ & 1.13451$\, 10^{-3}$ \\
  \midrule
  $\Sigma_{s,0}^{g\rightarrow 1}$ & 6.91158$\, 10^{-2}$ & 0. & 0. & 0.\\
  $\Sigma_{s,0}^{g\rightarrow 2}$ & 4.04132$\, 10^{-2}$ & 2.30626$\, 10^{-1}$ & 0. & 0.\\
  $\Sigma_{s,0}^{g\rightarrow 3}$ & 2.68621$\, 10^{-3}$ & 9.57027$\, 10^{-3}$ & 3.48414$\, 10^{-1}$ & 0.\\
  $\Sigma_{s,0}^{g\rightarrow 4}$ & 0. & 1.99571$\, 10^{-7}$ & 1.27195$\, 10^{-3}$ & 3.63631$\, 10^{-1}$\\
  
  \bottomrule[\heavyrulewidth] 
  \end{tabular}
    \caption{Cross sections in the Radial and Inner Blanket regions} 
  \label{tab:data-radial-blanket}
\end{table}

\begin{table}[htbp]
  \centering 
  \begin{tabular}{cccccc} 
  \toprule[\heavyrulewidth]
  \toprule[\heavyrulewidth]
   & \multicolumn{4}{c}{Radial Reflector}  \\
   \cmidrule(r){2-5} 
  \textbf{Group} & \textbf{1} & \textbf{2} & \textbf{3} &  \textbf{4}  \\ 
  \midrule
  $\Sigma_t$        & 1.71748$\, 10^{-1}$ & 2.17826$\, 10^{-1}$ & 4.47761$\, 10^{-1}$ & 7.95199$\, 10^{-1}$ \\
  $\nu\Sigma_f$     & 0. & 0. & 0. & 0. \\
  \midrule
  $\Sigma_{s,0}^{g\rightarrow 1}$ & 1.23352$\, 10^{-1}$ & 0. & 0. & 0.\\
  $\Sigma_{s,0}^{g\rightarrow 2}$ & 4.61307$\, 10^{-2}$ & 2.11064$\, 10^{-1}$ & 0. & 0.\\
  $\Sigma_{s,0}^{g\rightarrow 3}$ & 1.13217$\, 10^{-3}$ & 6.27100$\, 10^{-3}$ & 4.43045$\, 10^{-1}$ & 0.\\
  $\Sigma_{s,0}^{g\rightarrow 4}$ & 0. & 1.03831$\, 10^{-6}$ & 2.77126$\, 10^{-3}$ & 7.89497$\, 10^{-1}$\\
  
  \bottomrule[\heavyrulewidth] 
  \end{tabular}
    \caption{Cross sections in the Radial Reflector region } 
  \label{tab:data-radial-reflector}
\end{table}

\begin{table}[htbp]
  \centering 
  \begin{tabular}{cccccc} 
  \toprule[\heavyrulewidth]
  \toprule[\heavyrulewidth]
   & \multicolumn{4}{c}{Axial Blanket}  \\
   \cmidrule(r){2-5} 
  \textbf{Group} & \textbf{1} & \textbf{2} & \textbf{3} &  \textbf{4}  \\ 
  \midrule
  $\Sigma_t$        & 1.16493$\, 10^{-1}$ &  2.20521$\, 10^{-1}$ &  3.44544$\, 10^{-1}$ &  3.88356$\, 10^{-1}$ \\
  $\nu\Sigma_f$     & 1.31770$\, 10^{-2}$ &  1.26026$\, 10^{-4}$ &  1.52380$\, 10^{-4}$ &  7.87302$\, 10^{-4}$ \\
  \midrule
  $\Sigma_{s,0}^{g\rightarrow 1}$ & 7.16044$\, 10^{-2}$ &  0. &  0. &  0. \\
  $\Sigma_{s,0}^{g\rightarrow 2}$ & 3.73170$\, 10^{-2}$ &  2.10436$\, 10^{-1}$ &  0. &  0.\\
  $\Sigma_{s,0}^{g\rightarrow 3}$ & 2.21707$\, 10^{-3}$ &  8.59855$\, 10^{-3}$ &  3.37506$\, 10^{-1}$ &  0. \\
  $\Sigma_{s,0}^{g\rightarrow 4}$ & 0. &  6.68299$\, 10^{-7}$ &  1.68530$\, 10^{-3}$ &  3.74886$\, 10^{-1}$ \\
  
  \bottomrule[\heavyrulewidth] 
  \end{tabular}
    \caption{Cross sections in the Axial Blanket region} 
  \label{tab:data-axial-blanket}
\end{table}

\begin{table}[htbp]
  \centering 
  \begin{tabular}{cccccc} 
  \toprule[\heavyrulewidth]
  \toprule[\heavyrulewidth]
   & \multicolumn{4}{c}{Axial Reflector}  \\
   \cmidrule(r){2-5} 
  \textbf{Group} & \textbf{1} & \textbf{2} & \textbf{3} &  \textbf{4}  \\ 
  \midrule
  $\Sigma_t$        & 1.65612$\, 10^{-1}$ &  1.66866$\, 10^{-1}$ &  2.68633$\, 10^{-1}$ &  8.34911$\, 10^{-1}$ \\
  $\nu\Sigma_f$     & 0. &  0. &  0. &  0.\\
  \midrule
  $\Sigma_{s,0}^{g\rightarrow 1}$ & 1.15653$\, 10^{-1}$ &  0. &  0. &  0.\\
  $\Sigma_{s,0}^{g\rightarrow 2}$ & 4.84731$\, 10^{-2}$ &  1.60818$\, 10^{-1}$ &  0. &  0.\\
  $\Sigma_{s,0}^{g\rightarrow 3}$ & 8.46495$\, 10^{-4}$ &  5.64096$\, 10^{-3}$ &  2.65011$\, 10^{-1}$ &  0.\\
  $\Sigma_{s,0}^{g\rightarrow 4}$ & 0. &  6.57573$\, 10^{-7}$ &  2.41755$\, 10^{-3}$ &  8.30547$\, 10^{-1}$\\
  
  \bottomrule[\heavyrulewidth] 
  \end{tabular}
    \caption{Cross sections in the Axial Reflector region} 
  \label{tab:data-axial-relector}
\end{table}

\begin{table}[htbp]
  \centering 
  \begin{tabular}{cccccc} 
  \toprule[\heavyrulewidth]
  \toprule[\heavyrulewidth]
   & \multicolumn{4}{c}{Empty Matrix}  \\
   \cmidrule(r){2-5} 
  \textbf{Group} & \textbf{1} & \textbf{2} & \textbf{3} &  \textbf{4}  \\ 
  \midrule
  $\Sigma_t$        & 1.36985$\, 10^{-2}$ &  1.69037$\, 10^{-2}$ &  3.12271$\, 10^{-2}$ &  6.29537$\, 10^{-2}$ \\
  $\nu\Sigma_f$     & 0. &  0. &  0. &  0. \\
  \midrule
  $\Sigma_{s,0}^{g\rightarrow 1}$ & 9.57999$\, 10^{-3}$ &  0. &  0. &  0.\\
  $\Sigma_{s,0}^{g\rightarrow 2}$ & 3.95552$\, 10^{-3}$ &  1.64740$\, 10^{-2}$ &  0. &  0.\\
  $\Sigma_{s,0}^{g\rightarrow 3}$ & 8.80428$\, 10^{-3}$ &  3.91394$\, 10^{-4}$ &  3.09104$\, 10^{-2}$ &  0.\\
  $\Sigma_{s,0}^{g\rightarrow 4}$ & 0. &  7.72254$\, 10^{-8}$ &  1.77293$\, 10^{-4}$ &  6.24581$\, 10^{-2}$\\
  
  \bottomrule[\heavyrulewidth] 
  \end{tabular}
    \caption{Cross sections in the Empty Matrix region } 
  \label{tab:data-empty-matrix}
\end{table}

\begin{table} [htbp]
  \centering 
  \begin{tabular}{cccccc} 
  \toprule[\heavyrulewidth]
  \toprule[\heavyrulewidth]
   & \multicolumn{4}{c}{Control Rod}  \\
   \cmidrule(r){2-5} 
  \textbf{Group} & \textbf{1} & \textbf{2} & \textbf{3} &  \textbf{4}  \\ 
  \midrule
  $\Sigma_t$        & 1.84333$\, 10^{-1}$ &   3.66121$\, 10^{-1}$ &   6.15527$\, 10^{-1}$ &   1.09486 \\
  $\nu\Sigma_f$     & 0. &   0. &   0. &   0. \\
  \midrule
  $\Sigma_{s,0}^{g\rightarrow 1}$ & 1.34373$\, 10^{-1}$ &   0. &   0. &   0. \\
  $\Sigma_{s,0}^{g\rightarrow 2}$ & 4.37775$\, 10^{-2}$ &   3.18582$\, 10^{-1}$ &   0. &   0. \\
  $\Sigma_{s,0}^{g\rightarrow 3}$ & 2.06054$\, 10^{-4}$ &   2.98432$\, 10^{-2}$ &   5.19591$\, 10^{-1}$ &   0.\\
  $\Sigma_{s,0}^{g\rightarrow 4}$ & 0. &   8.71188$\, 10^{-7}$ &   7.66209$\, 10^{-3}$ &   6.18265$\, 10^{-1}$\\
  
  \bottomrule[\heavyrulewidth] 
  \end{tabular}
    \caption{Cross sections in the (empty) Control Rod regions} 
  \label{tab:data-CR}
\end{table}

\begin{table}[htbp] 
  \centering 
  \begin{tabular}{cccccc} 
  \toprule[\heavyrulewidth]
  \toprule[\heavyrulewidth]
   & \multicolumn{4}{c}{Na Filled CRP}  \\
   \cmidrule(r){2-5} 
  \textbf{Group} & \textbf{1} & \textbf{2} & \textbf{3} &  \textbf{4}  \\ 
  \midrule
  $\Sigma_t$        & 6.58979$\, 10^{-2}$ &   1.09810$\, 10^{-1}$ &   1.86765$\, 10^{-1}$ &   2.09933$\, 10^{-1}$ \\
  $\nu\Sigma_f$     & 0. &   0. &   0. &   0.\\
  \midrule
  $\Sigma_{s,0}^{g\rightarrow 1}$ & 4.74407$\, 10^{-2}$ &   0. &   0. &   0.\\
  $\Sigma_{s,0}^{g\rightarrow 2}$ & 1.76894$\, 10^{-2}$ &   1.06142$\, 10^{-1}$ &   0. &   0.\\
  $\Sigma_{s,0}^{g\rightarrow 3}$ & 4.57012$\, 10^{-4}$ &   3.55466$\, 10^{-3}$ &   1.85304$\, 10^{-1}$ &   0.\\
  $\Sigma_{s,0}^{g\rightarrow 4}$ & 0. &   1.77599$\, 10^{-7}$ &   1.01280$\, 10^{-3}$ &   2.08858$\, 10^{-1}$\\
  
  \bottomrule[\heavyrulewidth] 
  \end{tabular}
    \caption{Cross sections in the sodium filled Control Rod regions}  
  \label{tab:data-CRP}
\end{table}

\newpage

\bibliographystyle{siam}
\bibliography{manuscript}

\end{document}